\documentclass[3p,times]{elsarticle_preprint}

\usepackage{amssymb}

\usepackage{amsmath} 
\usepackage{amsfonts}
\usepackage{textcomp}
\usepackage{mathrsfs}
\usepackage{float}
\usepackage{caption}
\usepackage{subfigure}
\usepackage{mathabx}
\usepackage{graphicx}
\usepackage{verbdef}
\usepackage{mathabx}

\usepackage{listings}
\usepackage{fancyvrb}
\usepackage{caption}
\usepackage{enumitem}
\usepackage{cprotect}
\usepackage{verbatim}
\usepackage{amsmath}
\usepackage{tikz}

\journal{R. Cavoretto, A. De Rossi, E. Perracchione}

\begin{document}
\newcommand{\NN}{\mathbb{N}}
\newcommand{\ZZ}{\mathbb{Z}}
\newcommand{\QQ}{\mathbb{Q}}
\newcommand{\RR}{\mathbb{R}}
\newcommand{\CC}{\mathbb{C}}
\newcommand{\II}{\mathbb{I}}

\newcommand{\Cross}{\mathbin{\tikz [x=1.4ex,y=1.4ex,line width=.1ex] \draw (-0.5,-0.5) -- (1.5,1.5) (-0.5,1.5) -- (1.5,-0.5);}}%

\newtheorem{osse}{Remark}[section]
\newtheorem{prob}{Problem}[section]
\newtheorem{Def}{Definition}[section]
\newtheorem{theorem}{Theorem}[section]

\begin{frontmatter}



\title{Efficient computation of partition of unity interpolants through a block-based searching technique}


\author{R. Cavoretto\corref{cor1}}
\ead{roberto.cavoretto@unito.it}
\cortext[cor1]{Corresponding author.}

\author{A. De Rossi}
\ead{alessandra.derossi@unito.it}

\author{E. Perracchione}
\ead{emma.perracchione@unito.it}

\address{Department of Mathematics \lq\lq G. Peano\rq\rq, University of Torino, via Carlo Alberto 10, I--10123 Torino, Italy}

\begin{abstract}
In this paper we propose a new efficient interpolation tool, extremely suitable for large scattered data sets. The partition of unity method is used and  performed by blending Radial Basis Functions (RBFs) as local approximants and using locally supported weight functions. In particular we present a new space-partitioning data structure based on a partition of the underlying generic domain in blocks. This approach allows us to examine only a reduced number of blocks in the search process of the nearest neighbour points, leading to an optimized searching routine. Complexity analysis and numerical experiments in two- and three-dimensional interpolation support our findings. Some applications to geometric modelling are also considered. Moreover, the associated software package written in \textsc{Matlab} is here discussed and made available to the scientific community. 
\end{abstract}



\end{frontmatter}


\section{Introduction} \label{intro}

Meshfree methods are popular tools for solving problems of interpolation and numerical resolution of differential equations. They take advantage of being flexible with respect to geometry, easy to implement in higher dimensions, and can also provide high order convergence. Recently, in approximation theory a specific method has been proved to be effective for interpolation of large scattered data sets, the partition of unity method. Its origin can be found in the context of partial differential equations (PDEs) \cite{Babuska97,Melenk96}. In scattered data interpolation it is implemented using RBFs as local approximants, since this is the most efficient tool for interpolation of scattered data \cite{Fasshauer}.  The main disadvantage of radial kernel-based method is the computational cost associated with the solution of (usually) large linear systems, therefore recent researches have been directed towards a change of the basis, either rendering them more stable, or considering a local method involving RBFs (see e.g. \cite{Cavoretto15bb,DeMarchi15,Fasshauer12,Fornberg11,Pazouki11,Safdari}). Here we focus on the localized RBF-based partition of unity approximation. As the name of the partition of unity method suggests, in such local approach, the efficient organization of scattered data is the crucial step. Precisely, in literature, techniques as \emph{kd-trees}, which allow to partition data in a $k$-dimensional space, and related searching procedures have already been designed \cite{Arya98,Cavoretto15a,deBerg97,Fasshauer,Wendland05}. Even if such techniques enable us to work with high dimensions, they are not specifically implemented for the partition of unity method.

In this paper, starting from the results shown in \cite{Cavoretto12a,Cavoretto14a,Cavoretto14c}, where efficient searching procedures based on the partition of underlying domains in strips or crossed strips are considered, we propose a versatile software for bivariate and trivariate interpolation which makes use of a new partitioning structure, named \emph{block-based partitioning structure}, and a novel related searching procedure. It strictly depends on the size of the partition of unity subdomains. Such technique allows to deal with a truly large number of data with a relatively low computational complexity. 

More precisely, our procedure for bivariate  and trivariate interpolation consists in covering, at first, the reconstruction domain with several non-overlapping small squares or cubes, named \emph{blocks}. Then the usually large scattered data set is distributed among the different blocks by recursive calls to a sorting routine. Once the scattered data are stored in such blocks, an optimized searching procedure is performed enabling us to solve the local interpolation problems arising from the domain decomposition. 
Specifically, such structure, built ad hoc for the partition of unity method, enables us to run the searching procedure in constant time complexity, independently from the initial number of nodes. An extensive complexity analysis supports our findings and moreover comparisons with other common techniques, as kd-trees, will be carried out. Interpolating large scattered data sets using procedures competitive with the most advanced techniques is thus our main purpose.

A second meaningful feature of our procedures is the flexibility with respect to the problem geometry.
In general, in literature the scattered data interpolation problem is considered in very simple and regular domains, such as squares or cubes \cite{Fasshauer,Fasshauer_slide}. This approach is limiting in the context of meshfree methods because of  the versatility of the meshless technique with respect to domains having different shapes. Instead in this work, our aim is to provide an automatic software that allows to solve  scattered data interpolation  problems in generic domains. Specifically, here we focus on convex domains. This choice is due to the fact that our scope consists in solving interpolation problems in domains which are, in general, a priori unknown, i.e.  problems arising from applications \cite{C-D-P-V,Sabetta}.

In what follows, in order to point out the versatility of the software, we will investigate several applications of such algorithm. For 2D data sets we stress the importance of having such versatile tool in biomathematics, presenting a short sketch about the reconstruction of the attraction basins \cite{C-D-P-V}. The same approach can also be employed in the approximation of the so-called \emph{sensitivity surfaces} \cite{Sabetta}. Then, for 3D data sets, we analyze the problem of modeling implicit surfaces via partition of unity interpolation \cite{Ohtake,Wendland02}. It is known that the reconstruction of 3D objects is computationally expensive because of the large amount of data. Thus, the importance of having an efficient partitioning structure in such framework follows.

The paper is organized as follows. In Section \ref{prelim} we recall theoretical preliminaries on local RBF-based partition of unity approximation. In Section \ref{block_algs}, we describe in detail the block-based partition of unity algorithms for bivariate and trivariate interpolation, which are based on the use of the new block-based partitioning and searching procedures. Computational complexity of these interpolation algorithms is then analyzed in Section \ref{comp_cost}. In Section \ref{num_exp} we report numerical experiments devoted to point out the accuracy of our algorithms. Section \ref{appl} contains some applications in biomathematics and CAGD. Section \ref{conclusion} deals with conclusions and future work. 
We point out that the algorithms are made available to the scientific community in a downloadable free software package: 
	\begin{center}
		http://hdl.handle.net/2318/158790.
	\end{center}

	\section{Preliminaries} \label{prelim}
	In this section we briefly review the partition of unity approximation based on a localized use of RBF interpolants. This computational technique is meshfree and effectively works with large sets of scattered data points \cite{Fasshauer,Wendland05}. 
	
	\subsection{RBF interpolation} 
	Given a set ${ \cal X}_N= \{ \boldsymbol{x}_i \in \mathbb{R}^{M},i=1, \ldots, N \}$ of $N$ distinct \textsl{data points}, also called \textsl{data sites} or \textsl{nodes}, in a domain $ \Omega \subseteq \mathbb{R}^{M}$, and a corresponding set $ {\cal F}_N= \{ f_i = f(\boldsymbol{x}_i)  ,i=1, \ldots, N \}$ of \textsl{data values} or \textsl{function values} obtained by possibly sampling any (unknown) function  $f: \Omega \longrightarrow \mathbb{R}$, the standard RBF interpolation problem consists in finding an interpolant $R:\Omega \longrightarrow \mathbb{R}$ of the form
	\begin{equation}
	R( \boldsymbol{x})= \sum_{i=1}^{N} c_i \phi (||  \boldsymbol{x} - \boldsymbol{x}_i  ||_2), \quad  \boldsymbol{x} \in \Omega,
	\label{rad1}
	\end{equation}
	where $||\cdot||_2$ is the Euclidean norm, and $ \phi: [0, \infty) \longrightarrow \mathbb{R}$ is a RBF \cite{Buhmann03,Iske11}. The coefficients $ \{ c_i \}_{i=1}^{N} $ are determined by enforcing the interpolation conditions
	\begin{equation}
	R( \boldsymbol{x}_i)=f_i, \quad i=1, \ldots, N.
	\label{int1}
	\end{equation}
	Imposing the conditions \eqref{int1} leads to a symmetric linear system of equations
	\begin{equation}
	\Phi \boldsymbol{c}= \boldsymbol{f},
	\label{sys1}
	\end{equation}
	where $\Phi_{ki} = \phi (||  \boldsymbol{x}_k - \boldsymbol{x}_i  ||_2)$, $k,i=1, \ldots, N$, $\boldsymbol{c}= [c_1, \ldots, c_N]^T$, and $  \boldsymbol{f} =[f_1, \ldots , f_N]^T$. When $\boldsymbol{c}$ is found by solving the system \eqref{sys1}, we can evaluate the RBF interpolant at a point $\boldsymbol{x}$ as
	\begin{align*}
	R(\boldsymbol{x}) = \boldsymbol{\phi}^T(\boldsymbol{x}) \boldsymbol{c},
	\end{align*}
	where $\boldsymbol{\phi}^T(\boldsymbol{x}) = [\phi (||  \boldsymbol{x} - \boldsymbol{x}_1  ||_2),\ldots,\phi (||  \boldsymbol{x} - \boldsymbol{x}_N  ||_2)]$.
	
The interpolation problem is well-posed, i.e. a solution to the problem exists uniquely, if and only if the matrix $\Phi$ is nonsingular. A sufficient condition to have nonsingularity is that $\Phi$ is positive definite. 

	\subsection{Partition of unity approximation} \label{pu}
	
	Let $\Omega \subseteq \RR^M$ be an open and bound\-ed domain, and let $\{\Omega_j\}_{j=1}^{d}$ be an open and bounded covering of $\Omega$ satisfying some mild overlap condition among the subdomains $\Omega_j$,  i.e. the overlap among the subdomains must be sufficient so that each interior point $\boldsymbol{x} \in \Omega$ is located in the interior of at least one subdomain $\Omega_j$. The set $I(\boldsymbol{x}) = \{j : \boldsymbol{x} \in \Omega_j \}$, for $\boldsymbol{x} \in \Omega$, is uniformly bounded on $\Omega$, with $ \Omega  \subseteq \bigcup_{j=1}^{d} \Omega_j$. 
		
	Associated with the subdomains we choose partition of unity weight functions $W_j$, i.e. a family of compactly supported, nonnegative and continuous functions subordinate to the subdomain $\Omega_j$, such that $ \sum_{j = 1}^d W_j( \boldsymbol{x})=1$ on $\Omega$ and ${\rm supp}(W_j)  \subseteq \Omega_j$. The global approximant is thus constructed as follows
	\begin{equation} 
	{\cal I}( \boldsymbol{x})= \sum_{j=1}^{d} R_j( \boldsymbol{x} ) W_j ( \boldsymbol{x}), \quad \boldsymbol{x} \in \Omega,
	\label{intg}
	\end{equation}
	where $R_j$ defines a local RBF interpolant on each subdomain $\Omega_j$ and $W_j: \Omega_j \longrightarrow \mathbb{R}$ is a partition of unity weight function. 
	
	According to \cite{Wendland02}, we assume to have a $ k$-stable partition of unity, i.e a family of nonnegative functions $ \{ W_j \}_{j=1}^{d}$,  with $W_j \in C^k ( \mathbb{R}^M)$, such that:
	\begin{itemize}
		\item[i.] supp$ ( W_j ) \subseteq  \Omega_j $,
		\item[ii.] $ \sum_{j=1}^{d} W_j( \boldsymbol{x})=1$ on $ \Omega$,
		\item[iii.] $||D^{ \beta} W_j ||_{L^{ \infty} ( \Omega_j)} \leq \frac{C_{ \beta} }{ \delta_j^{ | \beta|}},$ $  \forall \beta \in \mathbb{N}^{M}: | \beta |	 \leq k,$
		where $ \delta_j$ is the diameter of $\Omega_j$ and $C_{ \beta} > 0$ is a constant.
	\end{itemize}
	As nonnegative functions $W_j \in C^k ( \mathbb{R}^M)$, we consider Shepard's weight, i.e.,
	\begin{align*}
	W_j(\boldsymbol{x}) = \frac{\varphi_j(\boldsymbol{x})}{\sum_{k\in I(\boldsymbol{x})} \varphi_k(\boldsymbol{x})}, \quad j=1,\ldots,d,
	\end{align*}
	$\varphi_j(\boldsymbol{x})$ being compactly supported functions with support on $\Omega_j$ such as Wendland's functions \cite{Wendland05}.

	\begin{osse}
		If the functions $R_j$, $j=1,\ldots,d$, satisfy the interpolation conditions $ R_j( \boldsymbol{x}_i )= f( \boldsymbol{x}_i)$ for each $\boldsymbol{x}_i \in \Omega_j$, then the global partition of unity approximant inherits the interpolation property of the local interpolants \cite{Fasshauer}, i.e. 
		\begin{align*}
		{\cal I}( \boldsymbol{x}_i ) = \sum_{j=1}^{d} R_j( \boldsymbol{x}_i ) W_j ( \boldsymbol{x}_i) = \sum_{j \in I(\boldsymbol{x}_i)} f( \boldsymbol{x}_i ) W_j ( \boldsymbol{x}_i) = f( \boldsymbol{x}_i).
		\end{align*}
	\end{osse}
	
	In order to be able to formulate error bounds, we need some further assumptions on regularity of $\Omega_j$ and define the \textsl{fill distance}
	\begin{align}
	h_{ {\cal X}_N, \Omega} =  \sup_{ \boldsymbol{x} \in \Omega} \min_{ \boldsymbol{x}_i  \in {\cal X}_N} || \boldsymbol{x} - \boldsymbol{x}_i||_2.
	\label{fd}
	\end{align}
	Specifically, we require that an open and bounded covering $ \{ \Omega_j \}_{j=1}^{d}$ is \textsl{regular} for $( \Omega, {\cal X}_N)$. This means to fulfill the following properties \cite{Wendland02a}:
	\begin{enumerate}
		\item[i.]  for each $ \boldsymbol{ x} \in \Omega$, the number of subdomains $ \Omega_j$ with $ \boldsymbol{x} \in \Omega_j$ is bounded by a global constant $C$;
		\item[ii.]  there exists a constant $C_r > 0$ and an angle $\theta \in (0,\pi/2)$ such that every subdomain $ \Omega_j$ satisfies an interior cone condition with angle $\theta$ and radius $r = C_r h_{ {\cal X}_N, \Omega}$;
		\item[iii.]  the local fill distances $ h_{ {\cal X}_{N_j},\Omega_j}$ are uniformly bounded by the global fill distance $h_{{\cal X}_N, \Omega}$, where ${\cal X}_{N_j}= {\cal X}_N \cap \Omega_j$.
	\end{enumerate}
	
	\begin{osse}
		The assumptions above lead to the requirement that the number of subdomains is proportional to the number of data \cite{Wendland05}.
		The first property ensures that \eqref{intg} is actually a sum over at most $C$ summands. Moreover, it is crucial for an efficient evaluation of the global approximant that only a constant number of local interpolants has to be evaluated. It follows that it should be possible to locate those $C$ indices in constant time. The second and third properties are significant for estimating errors of RBF interpolants.
		\label{osservaz_1}
	\end{osse}
	
	After defining the space $C_{ \nu}^{k}  ( \mathbb{R}^{M} ) $ of all functions $f \in C^k$ whose derivatives of order $ | \beta |=k $
	satisfy $ D^{ \beta} f ( \boldsymbol{x} ) = {\cal O} ( || \boldsymbol{x} ||_2^{ \nu} ) $ for $ || \boldsymbol{x} ||_2 \longrightarrow 0$, we consider the following convergence result \cite{Fasshauer, Wendland05}:
	\begin{theorem}
		Let $ \Omega \subseteq  \mathbb{R}^M$ be open and bounded and suppose that $  {\cal X}_N= \{ \boldsymbol{x}_i  ,i=1, \ldots, N \} \subseteq \Omega$. Let $ \phi \in C_{ \nu}^{k}  ( \mathbb{R}^{M} ) $ be a
		strictly   positive definite function. Let $ \{ \Omega_j \}_{j=1}^{d}$ be a regular covering for  $( \Omega,  {\cal X}_N)$  and let $ \{ W_j \}_{j=1}^{d}$ be $k$-stable for $ \{ \Omega_j  \}_{j=1}^{d}$. Then the error between $ f \in \mathscr{N}_{ \phi} ( \Omega)$, where $ \mathscr{N}_{ \phi} = \textrm{span} \{ \phi(||\boldsymbol{x}-\cdot||_2), \boldsymbol{x} \in \Omega\}, $ is the native space of  $ \phi $,   and its partition of unity interpolant \eqref{intg} can be bounded by:
		$$
		| D^{ \beta} f( \boldsymbol{x} ) -  D^{ \beta} {\cal I}( \boldsymbol{x} ) | \leq C^{'} h_{  {\cal X}_N, \Omega}^{\frac{ k+ \nu }{2} - | \beta |} |f|_{{\cal N}_{ \phi} ( \Omega )},
		$$
		for all $ \boldsymbol{x} \in \Omega $ and all $ | \beta | \leq k/2$, where $C^{'}$ is a constant independent of $\boldsymbol{x}$, $f$ and  $\phi$.
		\label{th1}
	\end{theorem}

	\begin{osse}
	If we compare the result reported in Theorem \ref{th1} with the global error estimates shown in  \cite{Wendland05}, we can see that the partition of unity interpolant preserves the local approximation order for the global fit.  Thus, the partition of unity approach enables us to decompose a large problem into many small ones and, at the same time, ensures that the accuracy obtained for the local fits is carried over to the global interpolant. 
	\label{osservaz_2}
	\end{osse}
	
	\begin{osse}
	From Theorem  \ref{th1}, we can note that the interpolation error decreases together with the fill distance. Anyway, consistently with the trade-off or uncertainty \emph{principle} \cite{Schaback}, a conflict between theoretical accuracy and numerical stability may occur. In fact, if a large number of interpolation nodes is involved, the local RBF systems may suffer from ill-conditioning.  The latter is linked to the order of the basis functions and to the node distribution. Therefore, the ill-conditioning grows if  the fill distance decreases. In such case, in order to avoid numerical problems, for high density of interpolation  points,	we can use low-order basis functions or Compactly Supported RBFs (CSRBFs) \cite{Fasshauer}. 	More recently, however, several approximation techniques	have been proposed to have a stable computation with flat RBFs \cite{Fasshauer15}.
	\label{osservaz_3}
	\end{osse}
	
	\section{Block-based interpolation algorithms} \label{block_algs}
	This section is devoted to the presentation of the partition of unity algorithms for bivariate and trivariate interpolation, which make use of the new block-based partitioning structure and related optimized searching procedure. They allow us to efficiently find  all the points belonging to a given subdomain $\Omega_j$, which as in \cite{Cavoretto14a,Cavoretto15a,Fasshauer,Safdari,Shcherbakov} consists of circular or spherical patches (depending on whether $M=2$ or $3$).
	
	Here, since our main target is the interpolation of large scattered data,  in the partition of unity scheme we compute the local interpolants by means of CSRBFs. However, as it will be pointed out, this approach turns out to be very flexible and different choices of local approximants, either globally or compactly supported, are allowed.
	
	Since we are going to describe in detail our \textsc{Matlab} routines, in Table \ref{mfun}  we first summarize  the functions of the proposed software.

	In what follows we will use a common notation for the \textsc{Matlab} routines listed in Table \ref{mfun}. As example, \verbdef\demo{BlockBas-}\demo \\
	\verbdef\demo{edMD_Structure.m}\demo \verbdef\demo{spazio} denotes both the routines \verbdef\demo{BlockBased2D_Structure.m}\demo \verbdef\demo{spazio} and \verbdef\demo{BlockBased3D_Structure.m}\demo \verbdef\demo{spazio}.
	
    Moreover, for easiness of the reader, the steps of the bivariate ($M=2$) and trivariate ($M=3$) partition of unity method,  which makes use of the block-based data structure and employs CSRBFs, are shown as pseudo-code in the \verbdef\demo{PUM_MD_CSRBF Algorithm}\demo. 
  
   		\begin{table}[ht]
   			\begin{center}
   				\begin{tabular}{cc} 
   					\hline\noalign{\smallskip}
   					\verbdef\demo{PUM_2D_CSRBF.m}\demo & scripts performing the partition  \\
   					\verbdef\demo{PUM_3D_CSRBF.m}\demo    & of unity using CSRBFs \\
   					\noalign{\smallskip}
   					\verbdef\demo{BlockBased2D_Structure.m}\demo &   scripts that store points into the \\
   					\verbdef\demo{BlockBased3D_Structure.m}\demo &   different neighbourhoods  \\
   					\noalign{\smallskip}
   					\verbdef\demo{BlockBased2D_ContainingQuery.m}\demo & scripts performing\\
   					\verbdef\demo{BlockBased3D_ContainingQuery.m}\demo &  the containing query procedure \\
   					\noalign{\smallskip}
   					\verbdef\demo{BlockBased2D_RangeSearch.m}\demo & scripts that perform the\\
   					\verbdef\demo{BlockBased3D_RangeSearch.m}\demo &  range search procedure\\                                                                                                                                                                                                                                         
   					\noalign{\smallskip}
   					\verbdef\demo{BlockBased2D_DistanceMatrix.m}\demo & scripts that form the distance matrix\\
   					\verbdef\demo{BlockBased3D_DistanceMatrix.m}\demo &  of two sets of points for  CSRBFs  \\                                                                                                                                      
   					\hline 
   				\end{tabular}
   			\end{center}
   			\caption{The \textsc{Matlab} codes for the block-based partition of unity algorithms.}
   			\label{mfun}
   		\end{table}
   		
   		  In order to construct a flexible procedure, at first, we need to focus on the problem geometry, i.e. we need a sort of data pre-processing, enabling us to consider scattered data sites arbitrarily distributed in a domain $\Omega \subseteq \mathbb{R}^M$, with $M=2$ or $3$.

   \subsection{The problem geometry}
   \label{pr_ge}
   
   In this subsection we refer to the {\fontfamily{pcr} \selectfont Step 1} of the \verbdef\demo{PUM_MD_CSRBF Algorithm}\demo.

    Since our aim is to construct an automatic algorithm for solving the interpolation problem of scattered data points arbitrarily distributed in a (a priori unknown) domain $\Omega$, the most appropriate way to act is to settle $\Omega$ as the convex hull 
    defined by the data set ${\cal X}_N= \{\boldsymbol{x}_i,i=1,\ldots,N\}$.
    This phase allows to approximate the interpolant on the minimal set containing points which can be automatically detected. Such strategy is neither limiting  nor restrictive in any sense, in fact  if the domain is supposed to be known, any generalization is possible and straightforward \cite{Heryudono,Safdari}. 

    After computing the convex hull, we need to define several auxiliary structures, meaningful to construct a robust partitioning data structure.
	Thus, we define a rectangular bounding ${\cal R}$ of the domain $\Omega$ as
	\begin{equation}
	{\cal R}= \prod_{m=1}^{M}   \left[ \min_{i=1, \ldots, N}  x_{im}, \max_{i=1, \ldots, N}  x_{im} \right].
	\label{rettangolo}
	\end{equation}		
	As evident from \eqref{rettangolo}, ${\cal R}$  consists of a rectangle or a rectangular prism, depending on whether $M=2$ or $3$.
	
\begin{figure}
		\begin{center}
			\includegraphics[height=.30\textheight]{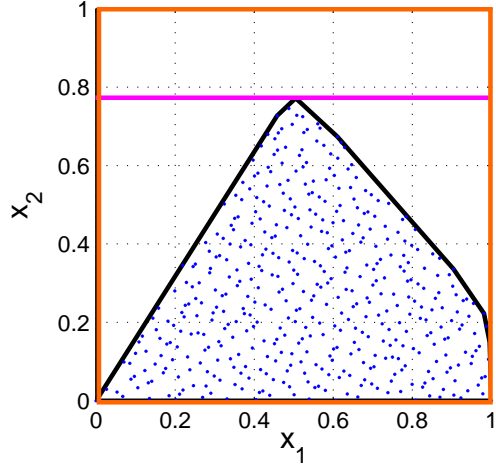} 
			\caption{An example of the problem geometry in a 2D framework: the set of data sites ${\cal X}_N$ (blue), 
				the convex hull $\Omega$ (black), the rectangle ${\cal R}$ containing  $\Omega$ (pink) and the  bounding box ${\cal L}$ (orange).} 
			\label{celle}
		\end{center}
	\end{figure}

	Moreover, in the problem geometry we consider a second auxiliary structure, known as \emph{bounding box}. This  is  a square or a cube (for $M=2$ or $3$ respectively) and is given by
	\begin{equation}
	{\cal L}=  \prod_{m=1}^{M}   \left[ \min_{m=1, \ldots, M} \left( \min_{i=1, \ldots, N}  x_{im} \right), \max_{m=1, \ldots, M} \left( \max_{i=1, \ldots, N}  x_{im} \right) \right].
	\end{equation}
	In order to fix the idea in a 2D framework refer to Figure \ref{celle}.

    \subsection{Definition of the partition of unity framework}
    \label{pu_def}
   	In this subsection we refer to the {\fontfamily{pcr} \selectfont Steps 2-5} of the \verbdef\demo{PUM_MD_CSRBF Algorithm}\demo.
    	
    The auxiliary structures previously defined, together with the convex hull, are useful to generate both the set of evaluation points ${\cal  E}_{s}= \{ \tilde{\boldsymbol{x}}_i   ,i=1, \ldots, s \} \subseteq \Omega$  and the set of partition of unity subdomain centres ${\cal C}_{d}= \{ \bar{\boldsymbol{x}}_i   ,j=1, \ldots, d \}  \subseteq \Omega$. 
    These sets are respectively obtained by generating $s_{\cal R}$  and $d_{\cal R}$ points as grids on ${\cal R}$. Then, they are automatically reduced by taking only those $s$ evaluation points and $d$ subdomain centres lying in $\Omega$.\VerbatimFootnotes \footnote{The points are automatically  reduced by the \verbdef\demo{inhull.m}\demo \verbdef\demo{spazio} function, provided by J. D'Errico, available at \cite{MCFE}.} 
    		
    As stated in Section \ref{prelim}, we require that the subdomains $\{ \Omega_j \}_{j=1}^{d}$ form an open, bounded and regular covering for $\Omega$. These assumptions affect the choice of the number of partition of unity centres and the one of the subdomain radius $\delta$.	
	Specifically, from Remark \ref{osservaz_1}, we know that the number of subdomains $d$ should be proportional to $N$. In particular, assuming to have a nearly uniform node distribution, $d$ is a suitable number of partition of unity subdomains on $ \Omega$ if \cite{Cavoretto15a,Fasshauer}
	\begin{equation}
	\dfrac{N}{d} \approx 2^M.
	\label{Nd}
	\end{equation} 
		
    Thus, denoting by $A_{\Omega}$  the area or the volume (for $M=2$ or $3$ respectively) of  the convex  hull\VerbatimFootnotes \footnote{The convex hull and its area or volume can be computed using the \textsc{Matlab}  routine \verbdef\demo{convhulln.m}\demo.}, from a simple proportion we find a suitable  number of subdomains initially generated on  ${\cal R}$
	\begin{align*}
	d_{{\cal R}}=\bigg \lfloor  \frac{\displaystyle 1}{\displaystyle 2} L \left(\frac{N}{A_{\Omega}}\right)^{1/M}\bigg \rfloor^M,
	\end{align*}
	where $L$ denotes the edge of the bounding box. So, the initial number of subdomains $d_{{\cal R}}$ is later reduced by taking only those  $d$ centres lying in $\Omega$ and, in this way, \eqref{Nd} is satisfied. 

	Also the subdomain radius $\delta$ must be carefully  chosen. In fact, the subdomains must be a covering  of the domain $\Omega$  satisfying the overlap condition (see Section \ref{prelim}).  
	The required property can be fulfilled taking as radius
	\begin{equation} 
	\delta = \frac{\displaystyle L \sqrt{2}}{\displaystyle \left(d_{{\cal R}}\right)^{1/M}}.
	\label{PU_radius}
	\end{equation}

	Moreover, we have to define the set of CSRBF centres ${\cal \hat{X}}_{\hat N}=\{ \hat{\boldsymbol{x}}_l, l=1, \ldots, \hat{N} \}$, which here, as in \cite{Fasshauer}, coincides with the set of data sites, i.e. ${\cal \hat{X}}_{\hat N} \equiv {\cal X}_{N}$.	

	 \subsection{The block-based partitioning structure}
	 In this subsection we refer to the \verbdef\demo{BlockBasedMD_Structure.m}\demo \verbdef\demo{spazio} routine (see {\fontfamily{pcr} \selectfont Steps 6-7} of the \verbdef\demo{PUM_MD_CSRBF }\demo\\
	 	\verbdef\demo{Algorithm}\demo).
	    
	Once the partition of unity subdomains are generated, the whole problem reduces to solve, for each subdomain, a local interpolation problem. Specifically, in the $j$-th local interpolation problem, only those data sites and evaluation points belonging to $\Omega_j$  are involved. Consequently, a  partitioning data structure and a related searching procedure must be employed to efficiently find the points located in each subdomain.
	
	In literature, to this scope the so-called kd-tree partitioning structures are commonly and widely used \cite{Arya98,deBerg97,Fasshauer}. A kd-tree, short for $k$-dimensional tree, is a space partitioning data structure for organizing points in a $k$-dimensional space. Here, since we have a $M$-dimensional space we should refer to such trees as Md-trees. But, in order to keep common notations we will go on calling them kd-trees. Following such approach, after building the tree structures for both data sites and evaluation points, the problem of  finding all points  belonging to a  given subdomain can be easily solved (see \cite{Cavoretto15a,Fasshauer,Wendland05} for details).
	
	In this work our aim is therefore to propose a new partitioning structure and, consequently, a new searching procedure built ad hoc for the interpolation purpose. The latter, besides being flexible as kd-tree, allows to find all the points belonging to a given subdomain $\Omega_j$ and turns out to be competitive in terms of computational time and cost. Such procedure is a partitioning data scheme based on  storing  points  into different  blocks, which are obtained  from the  subdivision of the bounding box auxiliary  structure ${\cal L}$ into several squares or cubes. 
	
	The number $q$ of blocks along one side of ${\cal L}$ is strictly  linked to the subdomain radius $\delta$ and is given by
	\begin{equation} 
	q= \bigg \lceil \frac{\displaystyle L}{\displaystyle \delta} \bigg \rceil.
	\label{q}
	\end{equation}
	From  \eqref{q}  we can deduce  that the block-based partitioning scheme depends on the construction of the partition of unity subdomains. In such framework, we will be able to get an efficient procedure to find the nearest points.
	
	Thus, after defining the width of the blocks as in \eqref{q},  
	we number blocks from $1$ to $q^M$. In a 2D context they are numbered from bottom to top, left to right, see Figure \ref{celle_fig}. For trivariate data sets, starting from the order shown in Figure \ref{celle_fig}, we continue numbering blocks along the quote as well.
		
    The block-based partitioning structure allows us to store both data sites and evaluation points in each of the $q^M$ blocks. At first, in such routine a sorting procedure is performed to order data sites along the first coordinate. Then recursive calls to the sorting routine are used to order data along the remaining coordinates enabling us to store points into the different blocks, i.e.:
	\begin{itemize}[leftmargin=*]
		\item[i.] the set $ {\cal X}_N$ is partitioned by the block-based partitioning structure into $q^M$ 
		subsets ${ \cal X}_{N_k}$, $k = 1, \ldots , q^M$, where
		${ \cal X}_{N_k}$  are the points stored  in the $k$-th block;
		\item[ii.] the set ${\cal E}_s$ is partitioned by the block-based partitioning structure  into $q^M$
		subsets ${ \cal  E}_{s_k}$, $k = 1, \ldots , q^M$,  where
		${ \cal E}_{s_k}$  are the points stored  in the $k$-th  block.
	\end{itemize}
	
	\begin{osse}
		In the  block-based partitioning structure, a sorting routine on the indices is needed. To this aim an optimized sorting procedure for integers is performed.\VerbatimFootnotes\footnote{The \textsc{Matlab} function
			\verbdef\demo{countingsort.m}\demo  \verbdef\demo{spazio} is a routine of the package called Sorting Methods, provided by
			B. Moore, available  at \cite{MCFE}.}
	\end{osse}
		 
	 		 \subsection{The block-based searching procedure}
	 		 In this subsection we refer to the \verbdef\demo{BlockBasedMD_RangeSearch.m}\demo \verbdef\demo{spazio} and 
	 		 \verbdef\demo{BlockBasedMD_ContainingQuery.m}\demo \verbdef\demo{spazio}  routines \verbdef\demo{spazio} (see {\fontfamily{pcr} \selectfont Step 8} of the \verbdef\demo{PUM_MD_CSRBF Algorithm}\demo).
	 		 
	 	After organizing in blocks data sites and evaluation points, in order to compute local fits, i.e. interpolants on each subdomain, we need to perform several procedures enabling us to answer the following queries, respectively known as \emph{containing query} and  \emph{range search}:
	 	\begin{itemize}[leftmargin=*]
	 		\item[i.] given a subdomain centre $ \bar{\boldsymbol{x}}_j \in \Omega$,
	 		return the $k$-th block containing $\bar{\boldsymbol{x}}_j$;
	 		\item[ii.] given a set of data points $\boldsymbol{x}_i \in {\cal X}_N$ 
	 		and a subdomain $\Omega_j$, find all points located in that 
	 		subdomain, i.e. $\boldsymbol{x}_i  \in {\cal X}_{N_j}={\cal X}_N \cap \Omega_j$.
	 	\end{itemize}
	 	
	 	Thus, we perform a containing query and a range search routines  based on the block-based partitioning scheme. To this aim,
	 	it is convenient to point out that in bivariate interpolation, blocks are generated by the intersection of two families of orthogonal strips. The former (numbered from $1$ to $q$) are parallel to the  $x_2$-axis, whereas the latter  (again numbered from $1$ to $q$) are parallel to the $x_1$-axis. 
	 	For 3D data sets blocks are generated by the intersection of three orthogonal rectangular prisms. In what follows, for simplicity, with abuse of notation we will continue to call such rectangular prisms with the term \lq\lq strips\rq\rq. Consistently with the bivariate case, the three families of strips are all numbered from $1$ to $q$. Moreover, the first family of strips is parallel to the $(x_2,x_3)$-plane, the second one is parallel to the $(x_1,x_3)$-plane and the last one is parallel to the $(x_1,x_2)$-plane.
	 	
\begin{figure}
			\begin{center}
			\includegraphics[height=.30\textheight]{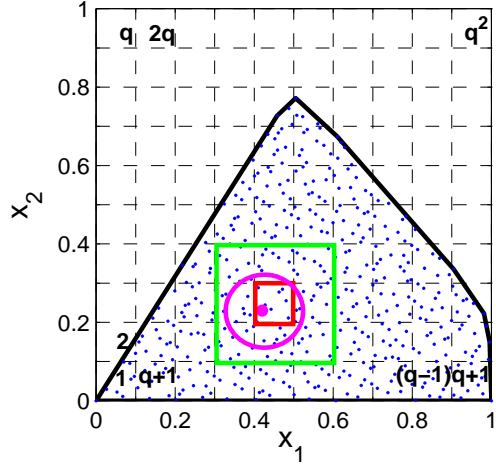} 
			\caption{An example of a 2D
			block-based  partitioning structure: the $k$-th block (red), a subdomain centre belonging to the $k$-th block (pink) and the neighbourhood set (green).}
			\label{celle_fig}
			\end{center}
			\end{figure}

	 	The block-based containing query, given a subdomain centre, returns the  index of the block containing  such centre. Thus, given a partition of unity centre $\boldsymbol{\bar{x}}_j$, if $k_m$ is the index of the strip parallel to the subspace of dimension $M-1$ generated by $x_p$, $p=1, \ldots,M$ and $p \neq m$, containing the $m$-th coordinate of $\boldsymbol{\bar{x}}_j$, then the index  of the  $k$-th block  containing the subdomain centre is
	 	\begin{align} 
	 	k=\sum_{m=1}^{M-1} \left( k_m-1 \right) q^{M-m}+k_M.
	 	\label{idx_2}
	 	\end{align}
	 	As example in a 2D framework, the subdomain centre plotted in Figure \ref{celle_fig} belongs to the $k$-th block, with $k=4q+3$; in fact here $k_1=5$ and $k_2=3$.  
	 		 		 	
	 	\begin{osse}
	 		The containing query routine is built to be a versatile tool usable for domains of any shape. In fact, this new structure turns out to be more flexible than the one proposed in \cite{Cavoretto14a,Cavoretto14c}, because it can be used for generic domains and not only for square or cube domains. 
	 	\end{osse}	
			
		After answering the first query, given a subdomain $\Omega_j$, the searching routine allows to:
		\begin{itemize}[leftmargin=*]
			\item[i.] find all data sites belonging to the subdomain $\Omega_j$;
			\item[ii.] determine all evaluation points belonging to the subdomain $\Omega_j$.
		\end{itemize}
		Specifically, supposing that the $j$-th subdomain centre belongs to the $k$-th block, the block-based searching procedure searches for all data lying in  the  $j$-th subdomain  among those lying in  the $k$-th neighbourhood, i.e.  in the $k$-th block and in its $3^M-1$ neighbouring blocks,  see Figure \ref{celle_fig}. 
		In particular, the partitioning structure based on blocks enables us to examine in the searching process at most $3^M-1$ blocks. In fact, when a block lies on the boundary of the bounding box, we reduce the number of neighbouring blocks to be considered.

	 		 \subsection{The computation of local block-based distance matrices and global interpolant}
	 		 In this subsection we refer to the \verbdef\demo{BlockBasedMD_DistanceMatrix.m}\demo  \verbdef\demo{spazio} routine (see {\fontfamily{pcr} \selectfont Steps 8a-10} of the \verbdef\demo{PUM_MD_CSRBF Algorithm}\demo).
	 		 			 	
	The  data sites and evaluation points belonging to the subdomain $\Omega_j$ are used to compute the local interpolation and evaluation matrices. In order to calculate the latter we have at first to compute the so-called \emph{distance matrices}.
	So, after that, the CSRBF is applied to the entire matrices obtaining the interpolation and evaluation matrices. In brief, referring to notation introduced in Section \ref{prelim}, this stage can be summarized as follows:
	\begin{enumerate}
		\item[1)] solving the local CSRBF linear system 
		\begin{align*}
		\Phi_j \boldsymbol{c}_j = \boldsymbol{f}_j,
		\end{align*}
		\item[2)] evaluating the local CSRBF interpolant
		\begin{align*}
		\boldsymbol{R}_j = \Phi_{j,eval} \boldsymbol{c}_j,
		\end{align*}
	\end{enumerate}
	where the index $j$ denotes the problem related to the $j$-th subdomain, while $\Phi_{j,eval}$ is the corresponding evaluation matrix \cite{Fasshauer}.
	
	Since we focus on CSRBFs, by properly scaling the support of the function, the local interpolation systems become sparse. Thus again, the block-based partitioning structure is used to efficiently find, for each  CSRBF centre, all data sites and evaluation points located within its support. As a consequence we compute  only few entries of the distance matrices. Finally, the local fits are accumulated into the global interpolant \eqref{intg}.
	
	On the opposite, in case of globally supported RBFs, since the entries of the distance matrices must be computed for each pair of points, building any partitioning structure is wasteful. Hence, in the \verbdef\demo{PUM_MD_CSRBF Algorithm}\demo, {\fontfamily{pcr} \selectfont Steps 8a-8b} should be skipped and the program \verbdef\demo{BlockBasedMD_DistanceMatrix.m}\demo \verbdef\demo{spazio} must be substituted by the function \verbdef\demo{DistanceMatrix.m}\demo, available in \cite{Fasshauer}. 

	
		\begin{table}
			\captionsetup{labelformat=empty}
			\begin{center}
				\begin{tabular}{p{16cm}*{1}{c}}
					\hline
					\vskip 0.01 cm 
					{\fontfamily{pcr} \selectfont INPUTS:} $N$, number of data; ${\cal X}_N=\{\boldsymbol{x}_i, i=1,\ldots,N\}$, set of data points; ${\cal F}_N=\{f_i, i=1,\ldots,N\}$, set of data 
					\vskip 0.08 cm 
					\hskip 2.2 cm values; $d_{{\cal R}}$, number of subdomains in ${\cal R}$; $s_{{\cal R}}$, number of evaluation points in ${\cal R}$.
					\vskip 0.12 cm 
					{\fontfamily{pcr} \selectfont OUTPUTS:} ${\cal A}_s=\{{\cal I}(\tilde{\boldsymbol{x}}_i), i=1,\ldots,s\}$, set of approximated values.
					\vskip 0.12 cm
					{\fontfamily{pcr} \selectfont Step 1:}
					Define the problem geometry, i.e. define $\Omega$ as the convex hull containing data sites, compute ${\cal R}$ and ${\cal L}$.
					\vskip 0.12 cm
					{\fontfamily{pcr} \selectfont Step 2:}
					A set ${\cal C}_d=\{\bar{\boldsymbol{x}}_j, j=1,\ldots,d\} \subseteq \Omega$ of subdomain points is constructed.
					\vskip 0.12 cm
					{\fontfamily{pcr} \selectfont Step 3:}
					A set ${\cal E}_s=\{\tilde{\boldsymbol{x}}_i, i=1,\ldots,s\} \subseteq \Omega$ of evaluation points is generated.
					\vskip 0.12 cm
					{\fontfamily{pcr} \selectfont Step 4:} Define the set of CSRBF centres ${\cal \hat{X}}_{\hat N}=\{\hat{\boldsymbol{x}}_l, l=1, \ldots, \hat{N} \} \subseteq \Omega$. Here  ${\cal \hat{X}}_{\hat N} \equiv  {\cal {X}}_N$.
					\vskip 0.12 cm
					{\fontfamily{pcr} \selectfont Step 5:} For each subdomain point $\bar{\boldsymbol{x}}_j$, $j=1,\ldots,d$, a  subdomain, whose radius  is given by \eqref{PU_radius}, is constructed.
					\vskip 0.12 cm
					{\fontfamily{pcr} \selectfont Step 6:} Compute the number $q$ of blocks  as in \eqref{q}.
					\vskip 0.12 cm
					{\fontfamily{pcr} \selectfont Step 7:}  The block-based data structures are built for the set ${\cal X}_N$ of data points and the set ${\cal E}_s$ of evaluation
					\vskip 0.08 cm
					\hskip 2.2 cm  points by using the  routine \verbdef\demo{BlockBasedMD_Structure.m}\demo \verbdef\demo|space|.   
					\vskip 0.12 cm
					{\fontfamily{pcr} \selectfont Step 8:}  For each subdomain $\Omega_j$, $j=1,\ldots,d$, the    \verbdef\demo{BlockBasedMD_ContainingQuery.m}\demo \verbdef\demo|space|  and the
					\vskip 0.08 cm
					\hskip 2.2 cm \verbdef\demo{BlockBasedMD_RangeSearch.m}\demo \verbdef\demo|space|  routines are performed allowing to:
					\vskip 0.08 cm
					\hskip 2.2 cm i. find all data points ${\cal X}_{N_j}$  belonging to the subdomain $\Omega_j$,
					\vskip 0.08 cm
					\hskip 2.2 cm ii. find all evaluation points  ${\cal E}_{s_j}$ belonging to the subdomain $\Omega_j$.
					\vskip 0.12 cm
					\hskip 1.6 cm {\fontfamily{pcr} \selectfont Step 8a:}  The block-based data structures are built for the set ${\cal X}_{N_j}$ of data points and the set ${\cal E}_{s_j}$ 
					\vskip 0.08 cm
					\hskip 4.2 cm  of evaluation points,  by using the \verbdef\demo{BlockBasedMD_Structure.m}\demo \verbdef\demo|space| routine.
					\vskip 0.12 cm
					\hskip 1.6 cm {\fontfamily{pcr} \selectfont Step 8b:} For each centre (of the basis function) $\boldsymbol{x}_i \in {\cal X}_{N_j}$, $i=1, \ldots, N_j$,  the
					\vskip 0.08 cm
					\hskip 4.2 cm  \verbdef\demo{BlockBasedMD_ContainingQuery.m}\demo  \verbdef\demo|space|  and the \verbdef\demo{BlockBasedMD_RangeSearch.m}\demo \verbdef\demo|space| 
					\vskip 0.08 cm
					\hskip 4.2 cm routines  are performed allowing to find:
					\vskip 0.08 cm
					\hskip 4.2 cm i. all data points ${\cal X}_{N_{j_i}}$ and
					\vskip 0.08 cm
					\hskip 4.2 cm ii. all evaluation points  ${\cal E}_{s_{j_i}}$  belonging to the support of the CSRBF centered at $\boldsymbol{x}_i$.
					\vskip 0.12 cm
					\hskip 1.6 cm {\fontfamily{pcr} \selectfont Step 9:} 
					\verbdef\demo{BlockBasedMD_DistanceMatrix.m}\demo \verbdef\demo|space| computes the interpolation and evaluation matrices
					\vskip 0.08 cm
					\hskip 4.2 cm and a local radial basis interpolant is formed.
					\vskip 0.12 cm
					{\fontfamily{pcr} \selectfont Step 10:} The local fits are accumulated into the global interpolant \eqref{intg}.
					\\[\smallskipamount]
					\hline
				\end{tabular}
			\end{center}
			\cprotect\caption{The \verbdef\demo{PUM_MD_CSRBF Algorithm}\demo. Routine performing the partition of unity method, using CSRBFs and employing the block-based partitioning structure and the related searching procedure.}
			\label{PUM_code}
		\end{table}
		
	\section{Complexity analysis} \label{comp_cost}
	
	In this section we point out the efficiency of our partitioning scheme. It will be proved that for  bivariate interpolation storing data sites and evaluation points requires
	$ {\cal O} (3/2 N \log N)$ and $ {\cal O} (3/2 s \log s)$  time complexity, respectively. While for 3D data sets
	the running times are $ {\cal O} (2 N \log N)$ and $ {\cal O} (2 s \log s)$ for storing data sites and evaluation points, respectively. 
	Moreover, when points are organized in blocks, for both 2D and 3D data sets the searching procedure can be computed in ${\cal O}(1)$ time complexity. This allows to perform a searching routine in a constant time, independently from the initial number of points.	A comparison with  kd-trees will be carried out (see Table \ref{tabella_costi}).

	\subsection{The block-based partitioning structure}
	The first part of the algorithm for partition of unity interpolation is a sort of data pre-processing which is not involved in complexity cost (see Subsections \ref{pr_ge} and \ref{pu_def}).
	
	Let us now focus on the partitioning structures used to organize the $N$  data sites in blocks. We remark that in the assessment of the total computational cost, analyzed in what follows in case of nodes, it must be added up the same  cost for storing the $s$ evaluation points.
	
	The partitioning structure employs the \emph{quicksort} routine which requires $ {\cal O} (n \log n)$ time complexity and $ {\cal O} ( \log n)$ space, where $n$ is the number of elements to be sorted. Specifically the block data structure is based on recursive calls to \verbdef\demo{sortrows.m}\demo \verbdef\demo{spazio}, which makes use of the \emph{quicksort} routine for sorting the nodes among the $M$ dimensions.

	To analyze the complexity of our procedures, we introduce the following notations, i.e.,
	\begin{displaymath}
	\begin{array}{ll}
	N_{k}^{(1)}: &   \textrm{number of data sites belonging to $q$ strips,}\\
	N_{k}^{(2)}: &   \textrm{number of data sites  belonging to $q^2$ strips.} 
	\end{array}
	\end{displaymath}
	
	Thus the computational cost  depending on the space dimension $M$ is:
	\begin{equation}
	{\cal O} \bigg( N \log{N} + \sum_{m=1}^{M-1} \sum_{k=1}^{q^m} N_k^{(m)} \log{N_k^{(m)}} \bigg).
	\label{CA1}
	\end{equation}
	Denoting by $N/q^{m}$ the average number of points lying in $q^m$ strips,  \eqref{CA1} can be estimated by
	\begin{equation}
	{\cal O} \bigg( N \log{N} + \sum_{m=1}^{M-1} N \log{ \frac{ \displaystyle N}{ \displaystyle q^m} } \bigg)
	\approx  {\cal O} \bigg( N \log{N} + \sum_{m=1}^{M-1} N \log{ \bigg( N \bigg(\frac{ \displaystyle \delta}{ \displaystyle L}\bigg)^m } \bigg)\bigg).
	\label{CA2}
	\end{equation}
	Now, from the definition of the partition of unity subdomains and neglecting the  
	constant terms, we obtain that \eqref{CA2} is approximately
	\begin{equation}
	{\cal O} \bigg( N \log{N} + \sum_{m=1}^{M-1} \frac{\displaystyle  M-m}{\displaystyle  M} N \log{ N} \bigg).
	\label{CA3}
	\end{equation}
	
		Moreover in the block-based partitioning scheme a sorting procedure on indices is employed to order them. Such routine  is performed with an optimized procedure for integers requiring  $ {\cal O} (n)$ time complexity,  where $n$ is the number of elements to be sorted. It follows that in the \lq\lq big O\rq\rq\ notation such cost turns out to be negligible in \eqref{CA3}.
	
	\subsection{The block-based searching procedure}
	To analyze the complexity of the 2D and 3D searching procedures, let $N_k$ be the number of data sites belonging to the $k$-th neighbourhood. Then, since  for each subdomain a quicksort procedure is used to order distances, the routine requires ${\cal O} ( N_k\log{N_k} )$ time complexity.  Observing that  the data sites  in a neighbourhood are about $N/(3q)^M$, the complexity can be estimated by
	\begin{equation}
	{\cal O} \bigg( \frac{ \displaystyle N}{\displaystyle (3q)^M} \log{ \frac{ \displaystyle N}{\displaystyle (3q)^M}}\bigg).
	\label{CA5}
	\end{equation}
	Taking into account the definitions of $q$ and $\delta$, \eqref{CA5} is approximately
	\begin{align}
	{\cal O} \bigg( \frac{ \displaystyle N 2^{M/2}}{\displaystyle 3^Md_{\cal R}} \log{ \frac{ \displaystyle N 2^{M/2}}{\displaystyle 3^Md_{\cal R}}}\bigg).
	\label{CA6}
	\end{align}
	Finally, substituting the definition of $d_{\cal R}$ in \eqref{CA6}, it is proved that
	\begin{equation}
	{\cal O} \bigg( A_{\Omega}\bigg(\frac{ \displaystyle 2\sqrt2}{\displaystyle 3 L}\bigg)^M \log{A_{\Omega}\bigg(\frac{ \displaystyle 2\sqrt2}{\displaystyle 3 L}\bigg)^M }\bigg) \approx {\cal O} (1).
	\label{CA7}
	\end{equation}
	
	The estimate \eqref{CA7} follows from the fact that we built a partitioning 
	structure strictly related to the size of the subdomains. 
	For this reason, in each partition of unity subdomain, the number of points
	is about constant, independently from the initial value $N$. 
	Thus using a number of blocks depending 
	both on the number and the size of such subdomains, the searching procedure involves a constant number of points, i.e. those belonging to a neighbourhood.
	\begin{osse}
		The same computational cost \eqref{CA3} and \eqref{CA7}, in case of CSRBFs, must be considered locally for each subdomain,  to build the sparse interpolation and evaluation matrices. In such steps we usually have a relatively small number of nodes $N_j$, with $N_j<< N$,  and evaluation points $s_j$, with $s_j<< s$, where the index $j$ identifies the $j$-th subdomain.
	\end{osse}
	
\begin{figure}
		\begin{center}
			\includegraphics[height=.26\textheight]{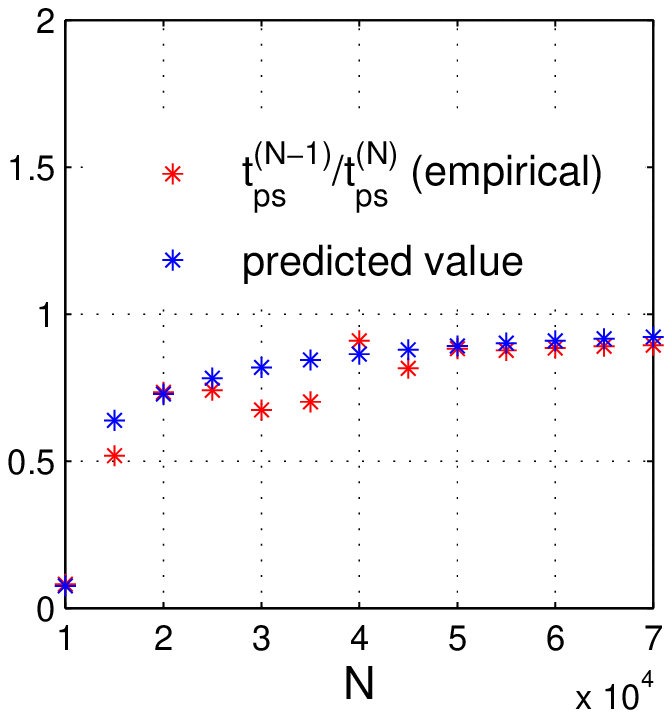} \hskip -.5cm 
			\includegraphics[height=.26\textheight]{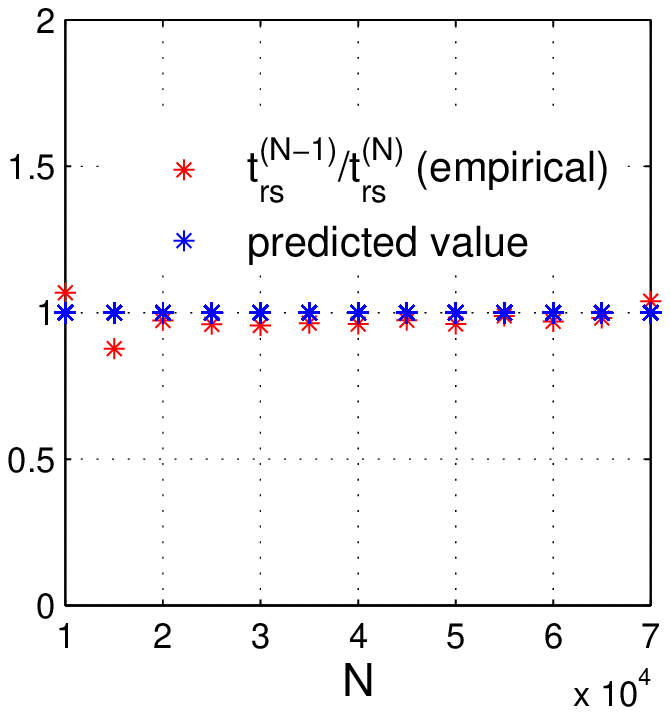} 
			\caption{
				The ratio between  empirical times  sampled at different consecutive 
				values of $N$ compared with the ratio of theoretical complexity costs, 
				for the block-based partitioning structure $t_{ps}$
				and for the  block-based searching procedure $t_{rs}$
				(left and right, respectively).
			}
			\label{time_comp}
		\end{center}
	\end{figure}

	All our findings are supported by numerical experiments shown in Figure \ref{time_comp}.  Here tests have been carried out on a Intel(R) Core(TM) i3 CPU M330 2.13 GHz processor.

	In Table \ref{tabella_costi}, we sum up the the total computational cost of the block-based partitioning and searching procedures, compared with kd-trees.
	
	\begin{table}[b!]
		\begin{center}
			\begin{tabular}{ccccc} 
				\hline\noalign{\smallskip}
				$M$	& Block-based &  kd-tree  & Block-based &  kd-tree \\
				& structure  & structure & search  & search   \\
				\noalign{\smallskip}
				\hline
				\noalign{\smallskip}
				2  & $ {\cal O} (3/2 N \log N)+$  & $ {\cal O} (2 N \log N)+$   & $ {\cal O} (1)$   &  ${\cal O} (\log N)+$   \\
				& $ {\cal O} (3/2 s \log s)$  &  $ {\cal O} (2 s \log s)$  &   &  ${\cal O} (\log s)$   \\
				\noalign{\smallskip}
				3 & $ {\cal O} (2 N \log N)+$   &  $ {\cal O} (3 N \log N)+$   & $ {\cal O} (1)$  & ${\cal O} (\log N)+$ \\	
				& $ {\cal O} (2 s \log s)$   & $ {\cal O} (3 s \log s)$   &     & ${\cal O} (\log s)$ \\
				\hline 
			\end{tabular}
		\end{center}
		\caption{Computational costs concerning block-based and the kd-tree routines.}
		\label{tabella_costi}
	\end{table}

	\subsection{Computation of local block-based distance matrices and global interpolant}
	Since the number of centres in each subdomain is bounded by a constant, we need ${ \cal O}(1)$ space and time for each subdomain to solve the local RBF interpolation  problems. In fact, to get the local interpolants, we have to solve $d$ linear systems of size $N_j \times N_j$, with $N_j << N$, thus requiring a constant running time ${\cal O}(N^3_j)$, $j =1, \ldots , d$, for each subdomain. Besides reporting the points in each subdomain in ${\cal O}(1)$, as the number $d$ of subdomains is bounded by ${\cal O}(N)$, this leads to ${\cal O}(N)$ space and time for solving all of them. Finally, we have to add up a constant number of local RBF interpolants to get the value of the global fit \eqref{intg}. This can be computed in ${\cal O}(1)$ time.

	\section{Numerical experiments} \label{num_exp}
	In our results we report errors obtained by running the algorithms on large scattered data sets located in convex hulls $\Omega \subseteq [0,1]^M$, for $M=2,3$. As interpolation points, we take uniformly random Halton data on the unit square or cube and then suitably reduced to $\Omega$.\VerbatimFootnotes\footnote{The Halton points are generated using the \textsc{Matlab} function \verbdef\demo{haltonseq.m}\demo \verbdef\demo{spazio},  provided by D. Dougherty, available at \cite{MCFE}.}   This choice  allows to make our tests repeatable.
	
Since we want to point out the efficiency of our partitioning routine, we also report  CPU times.
	
	\begin{osse}
To the best of our knowledge, the only fully available \textsc{Matlab} package for kd-trees, written by P. Vemulapalli, is given in \cite{MCFE}, but it is not optimally implemented. Thus a comparison on running times with our routines is not particularly meaningful. Moreover, for completeness,  we have to mention another package for kd-trees, written by  G. Shechter, available at \cite{MCFE}. It consists in dynamic libraries which are not executable in the recent  versions of  \textsc{Matlab} \cite{Fasshauer_slide,Kamranian}.
	\end{osse}
	
	To  point out the accuracy of our tests we will refer to 
	the \emph{maximum absolute error} (MAE) and the \emph{root mean square error}
	(RMSE), whose formulas are:
	\begin{align*} 
	MAE = \max_{1\leq i \leq s} |f(\tilde{\boldsymbol{x}}_i) - {\cal I}(\tilde{\boldsymbol{x}}_i)|, \qquad RMSE = \sqrt{\frac{1}{s}\sum_{i=1}^{s} |f(\tilde{\boldsymbol{x}}_i) - {\cal I}(\tilde{\boldsymbol{x}}_i)|^2}.
	\end{align*}
	
	Furthermore,  we also investigate two conditioning estimates, named the \emph{maximum conditioning number}
	and the \emph{average conditioning number}:
	\begin{align}
	\textrm{MaxCond}= \max_{1 \leq j \leq d} cond(\Phi_j), \qquad \textrm{AvCond}= \frac{1}{d} \sum_{j=1}^{d} cond(\Phi_j),
	\label{num_cond}
	\end{align}
	where $\Phi_j$ denotes the $j$-th matrix associated with the subdomain $\Omega_j$. More precisely, since the partition of unity method leads to solve  $d$ linear systems, to obtain a good conditioning
	estimate, in the right formula of \eqref{num_cond} we make an average among the conditioning numbers of the $d$ matrices.
	
	In these numerical experiments we focus on CSRBFs which might lead to sparse linear systems. So, for our propose we consider the compactly supported Wendland's $C^2$ function 
	\begin{equation}
	\phi(r)=(1-\varepsilon r)^{4}_{+} (4 \varepsilon r+1),
	\label{Wc2}
	\end{equation}
	where $ \varepsilon \in  \mathbb{R}^{+}$ is the shape parameter and $(\cdot)_{+}$ denotes the truncated power function. It follows that the  function \eqref{Wc2} is non negative for $r \in [0,1/\varepsilon ]$ and  strictly positive definite in $ \mathbb{R}^3$.

	\subsection{Results for bivariate interpolation}
	
\begin{figure}
		\centering
		\includegraphics[height=.26\textheight]{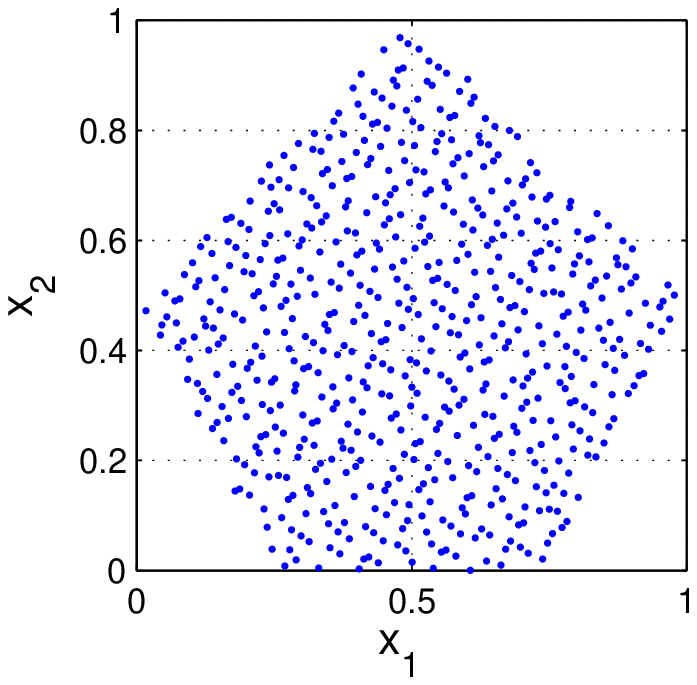} \hskip -0.5cm
		\includegraphics[height=.26\textheight]{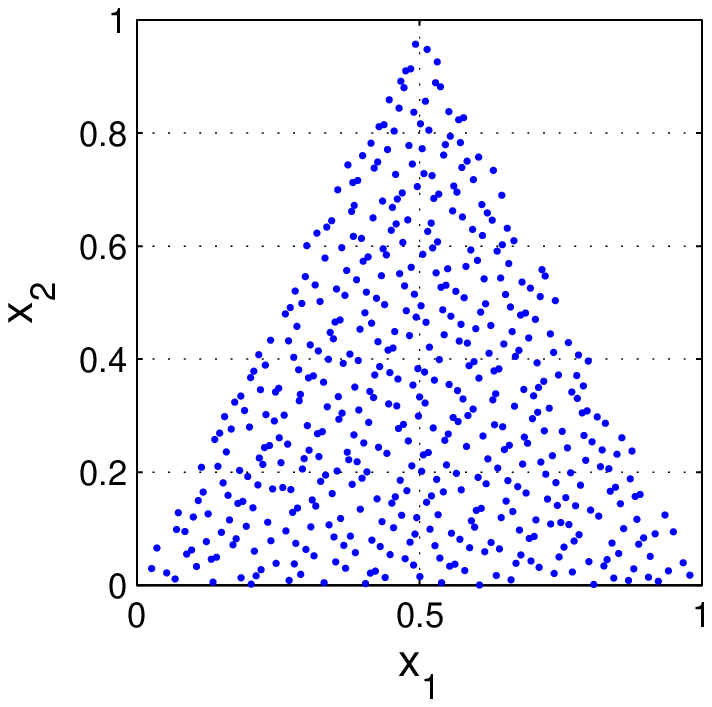}
		\caption{Examples of points in 2D convex hulls: pentagon with  $622$ nodes (left) and triangle with  $501$ nodes (right).}
		\label{fig_convex}
	\end{figure}

	In this subsection we focus on bivariate interpolation, analyzing performances of our algorithm and showing the numerical results obtained by considering five sets of Halton data points. These tests are carried out considering different convex domains, i.e. a triangle and a pentagon, see Figure \ref{fig_convex}.
	
	In the various experiments we investigate accuracy of the interpolation algorithm taking the data values by the well-known 2D Franke's  function $f_1$ and by the test function $f_2$:
	\begin{align*}
	f_1(x_1,x_2)&= \frac{3}{4}{\rm e}^{-\frac{(9x_1-2)^2+(9x_2-2)^2}{4}}+\frac{3}{4}{\rm e}^{-\frac{(9x_1+1)^2}{49}-\frac{9x_2+1}{10}} \\
	&+\frac{1}{2} {\rm e}^{-\frac{(9x_1-7)^2+(9x_2-3)^2}{4}}-\frac{1}{5} {\rm e}^{-(9x_1-4)^2-(9x_2-7)^2},
	\end{align*}
	\begin{align*}
	f_2(x_1,x_2)&= \frac{\displaystyle 1.25+ \cos{(5.4x_2)}}{6+6(3x_1-1)^2}.
	\end{align*}
	
	In Tables \ref{tabe_1} and \ref{tabe_3} we show the accuracy indicators of our algorithm considering several sets of  points for  pentagon and triangle, using $f_1$ and $f_2$ as test functions, respectively. These results are obtained taking  the shape parameter  $\varepsilon$ of \eqref{Wc2} equal to $0.5$ and a uniform grid of $40 \times 40$ evaluation points  on ${\cal R}$. Furthermore we also calculate the fill distance  \eqref{fd} and we estimate the empirical convergence rate via the formula:
	\begin{equation*}
	\lambda_k= \dfrac{\log(\textrm{RMSE}_{k-1}/\textrm{RMSE}_{k})}{\log(h_{{\cal X}_{N_{k-1}}}/h_{{\cal X}_{N_{k}}})}, \quad k=2,3, \ldots
	\end{equation*}
	where $\textrm{RMSE}_{k}$ is the error for the $k$-th numerical experiment, and $h_{{\cal X}_{N_{k}}}$ is the fill distance of the $k$-th computational mesh. Finally, in order to point out the efficiency, we also report in Tables \ref{tabe_1} and \ref{tabe_3} the CPU times (in seconds).

	\begin{table} 
		\begin{center}
			\begin{tabular}{cccccccc}
				\hline\noalign{\smallskip}
				$N$  & MAE & RMSE &  MaxCond & AvCond &  $ h_{{\cal X}_N}$ & $\lambda$ & $ t_{block}$\\
				\noalign{\smallskip}
				\hline
				\noalign{\smallskip}
				$622$  & $1.65{\rm E}-03$ & $1.40{\rm E}-04$ & $ 1.30 {\rm E}+07$ & $7.12 {\rm E}+06$ & $3.30{\rm E}-02$ & & $1.0$\\
				$2499$  & $5.02{\rm E}-04$ & $3.30{\rm E}-05$ & $ 1.72 {\rm E}+08$ & $4.82 {\rm E}+07$ & $1.76{\rm E}-02$& $2.29$ & $3.7$\\ 
				$9999$  & $4.33{\rm E}-05$ & $6.33{\rm E}-06$ & $ 1.92 {\rm E}+09$ & $5.46 {\rm E}+08$ & $9.12{\rm E}-03$ & $2.50$ & $9.1$\\ 
				$39991$  & $9.86{\rm E}-06$ & $1.25{\rm E}-06$ & $ 1.96 {\rm E}+10$ & $3.99 {\rm E}+09$ & $3.82{\rm E}-03$ & $1.86$ & $34.1$\\  
				$159994$   & $1.67{\rm E}-06$ & $3.05 {\rm E}-07$ & $ 1.74 {\rm E}+11$ & $3.56 {\rm E}+10$ & $2.02{\rm E}-03$ & $2.19$& $142.3$\\ 
				\hline
			\end{tabular}
		\end{center} 
		\caption{Errors, conditioning numbers, fill distances, convergence rates and CPU times (in seconds) for pentagon using  $f_1$ as test function and \eqref{Wc2} as local approximant with $\varepsilon = 0.5$.}
		\label{tabe_1}
	\end{table}

		\begin{table}[ht!]
		\begin{center}
			\begin{tabular}{cccccccc}
				\hline\noalign{\smallskip}
				$N$  & MAE & RMSE &  MaxCond & AvCond &  $ h_{{\cal X}_N}$ & $\lambda$ & $ t_{block}$\\
				\noalign{\smallskip}
				\hline
				\noalign{\smallskip}
				$501$  & $2.04{\rm E}-04$ & $2.60 {\rm E}-05$ & $ 1.24 {\rm E}+07$ & $6.98 {\rm E}+06$  & $3.35{\rm E}-02$& & $0.8$ \\
				$2008$  & $3.55{\rm E}-05$ & $5.41{\rm E}-06$ & $ 1.77 {\rm E}+08$ & $4.76 {\rm E}+07$ & $1.74{\rm E}-02$& $2.24$& $2.2$\\ 
				$7995$  & $2.32{\rm E}-05$ & $1.56{\rm E}-06$ & $ 1.77 {\rm E}+09$ & $5.47 {\rm E}+08$ & $9.16{\rm E}-03$ & $1.91$& $7.2$ \\ 
				$31999$  & $7.52{\rm E}-06$ & $4.59{\rm E}-07$ & $ 1.88 {\rm E}+10$ & $3.96 {\rm E}+09$ & $4.35{\rm E}-03$ & $1.63$ & $28.4$\\ 
				$128010$   & $5.75{\rm E}-07$ & $8.67{\rm E}-08$ & $ 1.67 {\rm E}+11$ & $3.55 {\rm E}+10$ & $2.15{\rm E}-03$& $2.32$& $114.3$\\ 
				\hline
			\end{tabular}
		\end{center} 
		\caption{Errors, conditioning numbers, fill distances, convergence rates and CPU times (in seconds) for triangle using  $f_2$ as test function and \eqref{Wc2} as local approximant with $\varepsilon = 0.5$.}
		\label{tabe_3}
	\end{table}

From the results shown in Tables \ref{tabe_1} and \ref{tabe_3}, we can see that, consistently with Remark \ref{osservaz_3}, the ill-conditioning grows in correspondence of a  decrease of the separation distance and of the errors. Furthermore, comparing the convergence rates reported  in Tables \ref{tabe_1} and \ref{tabe_3}, with the ones obtained for a global interpolant shown in \cite{Fasshauer}, we observe that the local convergence rates are carried over to the global interpolant. 
Hence, by means of the partition of unity method, together with the partitioning structure here proposed, we can efficiently and accurately decompose a large interpolation problem into many small ones (see Remark \ref{osservaz_2}). 
	
\begin{figure}
		\centering
		\includegraphics[height=.26\textheight]{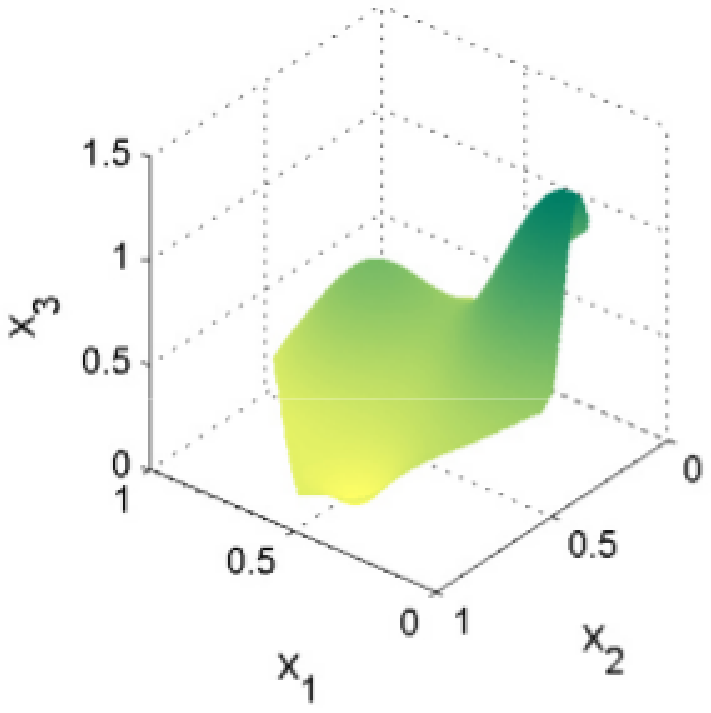} \hskip .5cm
		\includegraphics[height=.26\textheight]{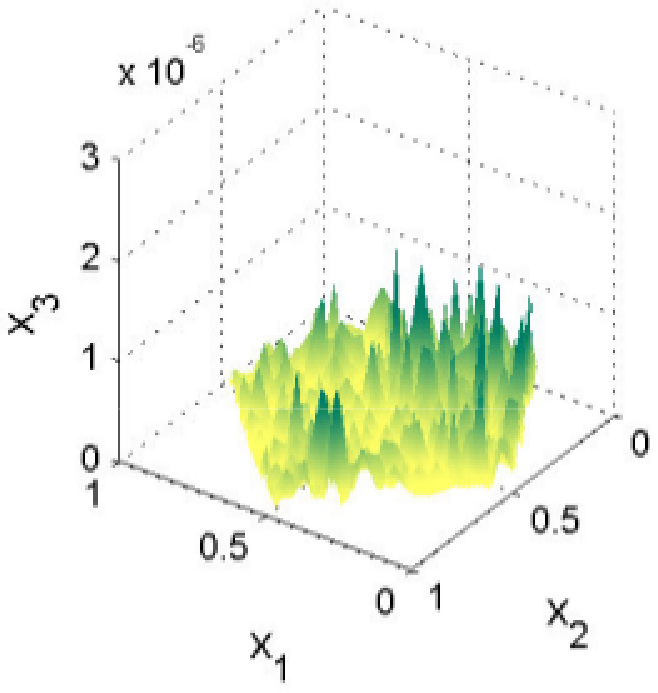}
		\includegraphics[height=.26\textheight]{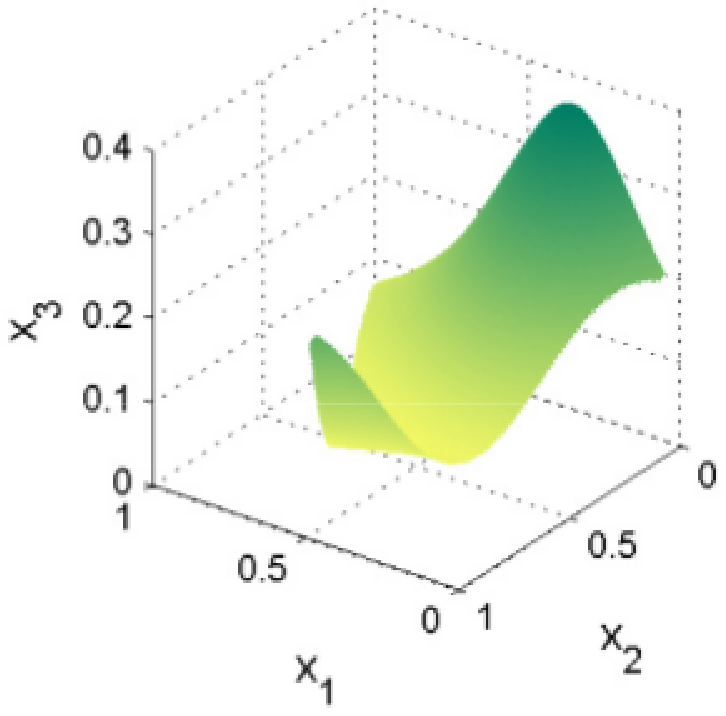} \hskip .5cm
		\includegraphics[height=.26\textheight]{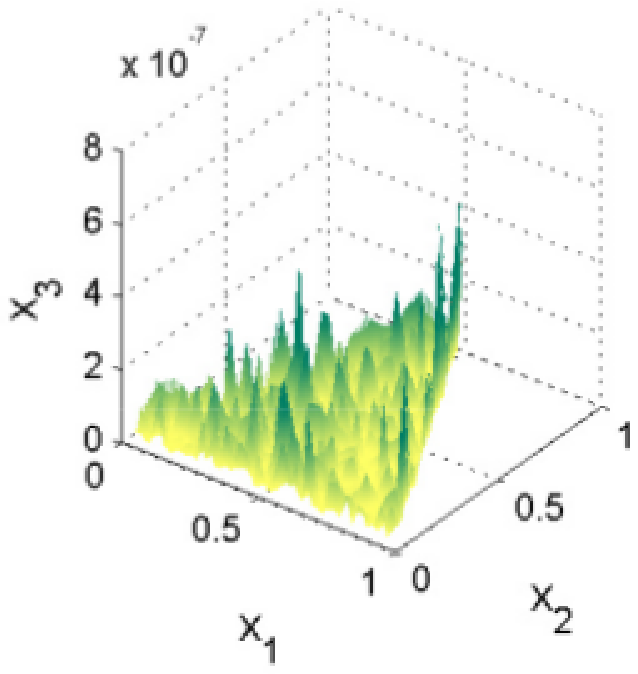} 
		\caption{The function $f_1$ (top left) and absolute errors (top right) defined on the pentagon with $159994$ nodes; the function $f_2$ (bottom left) and absolute errors (bottom right) defined on the triangle with $128010$ nodes.}
		\label{convex_franke}
	\end{figure}

	Moreover, in Figure \ref{convex_franke} we represent the two different test functions (left) and the absolute errors (right) computed on convex domains. 
		
	\subsection{Results for trivariate interpolation}
\begin{figure}
		\centering
		\includegraphics[height=.26\textheight]{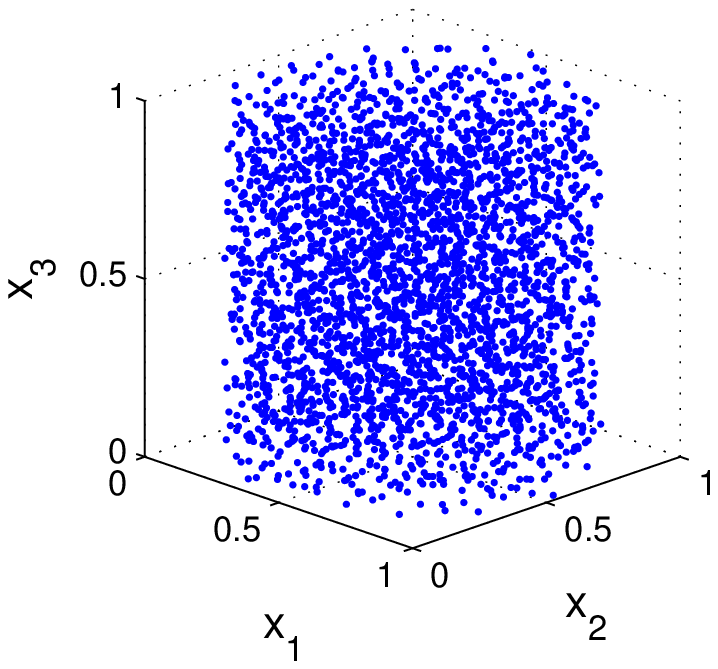} \hskip -.5cm
		\includegraphics[height=.26\textheight]{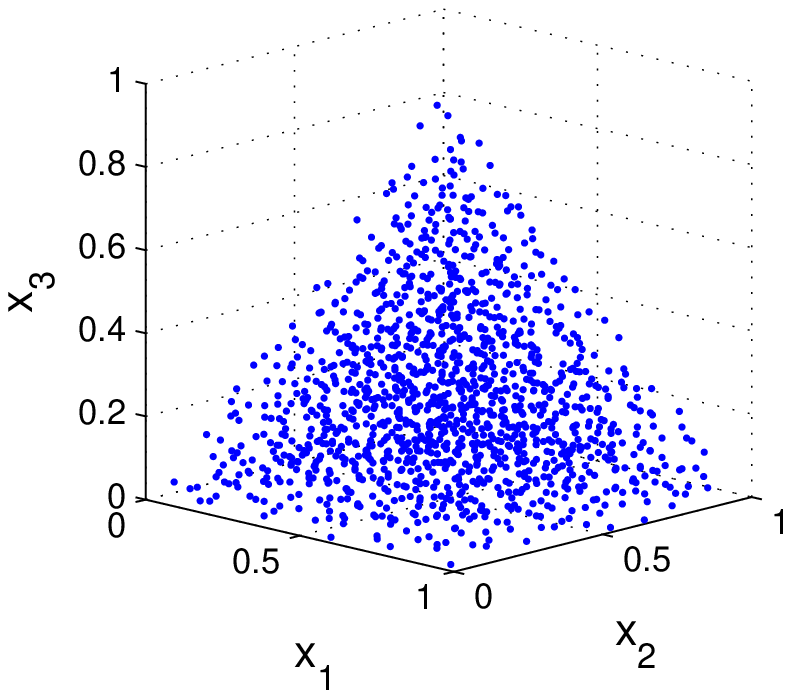}
		\caption{Examples of points in 3D convex hulls: cylinder with  $3134$ nodes (left) and pyramid with  $2998$ nodes (right).}
		\label{fig_convex2}
	\end{figure}

	In this subsection we instead report numerical results concerning trivariate interpolation. We analyze accuracy and efficiency 
	of the partition of unity algorithm for convex hulls, taking also in this case some sets of Halton scattered data points. Such points are located in a cylinder and in a pyramid, see Figure \ref{fig_convex2}.
	
	The trivariate test functions we consider in this subsection are the 3D Franke's function $f_3$ and the function $f_4$:
	\begin{eqnarray}
	f_3(x_1,x_2,x_3)&=& \frac{3}{4}{\rm e}^{-\frac{(9x_1-2)^2+(9x_2-2)^2+(9x_3-2)^2}{4}}+\frac{3}{4} {\rm e}^{-\frac{(9x_1+1)^2}{49}-\frac{9x_2+1}{10}-\frac{9x_3+1}{10}} \nonumber \\
	&+&\frac{1}{2} {\rm e}^{-\frac{(9x_1-7)^2+(9x_2-3)^2+(9x_3-5)^2}{4}}-\frac{1}{5} {\rm e}^{-(9x_1-4)^2-(9x_2-7)^2-(9x_3-5)^2}, \nonumber
	\end{eqnarray}
	\begin{eqnarray}
	f_4(x_1,x_2,x_3)&=& 4^3 x_1(1-x_1)x_2(1-x_2)x_3(1-x_3).  \nonumber
	\end{eqnarray}
	Tables \ref{tabe_5} and \ref{tabe_7} show the accuracy indicators of our algorithm considering several sets of Halton points for  cylinder and pyramid, using $f_3$ and $f_4$ as test functions, respectively. As earlier, these results are obtained taking  the shape parameter  $\varepsilon$ of \eqref{Wc2} equal to $0.5$.
	\begin{table}[ht!] 
		\begin{center}
			\begin{tabular}{cccccc}
				\hline\noalign{\smallskip}
				$N$  & MAE & RMSE &  MaxCond & AvCond \\
				\noalign{\smallskip}
				\hline
				\noalign{\smallskip}
				$3134$  & $5.94{\rm E}-03$ & $2.71{\rm E}-04$ & $ 7.65 {\rm E}+06$ & $3.74  {\rm E}+06$ \\ 
				$12551$  & $1.67{\rm E}-03$ & $6.00{\rm E}-05$ & $ 4.38{\rm E}+07$ & $2.07 {\rm E}+07$ \\ 
				$50184$  & $4.67{\rm E}-04$ & $2.27{\rm E}-05$ & $ 1.73 {\rm E}+08$ & $6.03 {\rm E}+07$ \\ 
				$200734$  & $1.22{\rm E}-04$ & $7.49{\rm E}-06$ & $9.86 {\rm E}+08$ & $3.35 {\rm E}+08$ \\ 
				$802865$   & $3.81{\rm E}-05$ & $2.91 {\rm E}-06$ & $ 3.10 {\rm E}+09$ & $1.07 {\rm E}+09$ \\ 
				\hline
			\end{tabular}
		\end{center} 
				\caption{Errors and conditioning numbers for cylinder using  $f_3$ as test function and \eqref{Wc2} as local approximant with $\varepsilon = 0.5$.}
		\label{tabe_5}
	\end{table}
	\begin{table} [ht!]
		\begin{center}
			\begin{tabular}{cccccc}
				\hline\noalign{\smallskip}
				$N$  & MAE & RMSE &  MaxCond & AvCond \\
				\noalign{\smallskip}
				\hline
				\noalign{\smallskip}
				$2998$  & $1.23{\rm E}-02$ & $6.40{\rm E}-04$ & $ 7.54 {\rm E}+06$ & $3.62  {\rm E}+06$ \\
				$12004$  & $3.38{\rm E}-03$ & $1.47{\rm E}-04$ & $ 4.49{\rm E}+07$ & $2.06 {\rm E}+07$ \\ 
				$47997$  & $7.33{\rm E}-04$ & $3.57{\rm E}-05$ & $ 1.73 {\rm E}+08$ & $6.03 {\rm E}+07$ \\ 
				$191981$  & $1.10{\rm E}-04$ & $9.93{\rm E}-06$ & $ 9.85 {\rm E}+08$ & $3.34 {\rm E}+08$ \\ 
				$767970$   & $4.63{\rm E}-05$ & $3.48 {\rm E}-06$ & $ 3.10 {\rm E}+09$ & $1.07 {\rm E}+09$ \\ 
				\hline
			\end{tabular}
		\end{center} 
				\caption{Errors and conditioning numbers for pyramid  using  $f_4$ as test function and \eqref{Wc2} as local approximant with $\varepsilon = 0.5$.}
		\label{tabe_7}
	\end{table}
	
	In Table \ref{tabe_6} we report the CPU times obtained by running  the block-based algorithm  for several sets of Halton data  in the  cylinder.  The results in  Table \ref{tabe_6} are obtained by considering a grid  of $20 \times 20 \times 20$ evaluation points  on ${\cal R}$. 
	Here, we omit the table concerning CPU times by varying $N$ for the pyramid because the behavior is similar to that outlined in Table \ref{tabe_6}.

	Also in the trivariate case, we register the pattern already discovered about ill-conditioning and accuracy. We can easily note that the ill-conditioning grows as the errors decrease. From this fact we can deduce that we have convergence. However, in this case we left out the computation of the convergence rates because the matrices become increasingly dense and computation requires lots of system memory.
	
	\begin{table}[ht!]
		\begin{center}
			\begin{tabular}{ccccc} 
				\hline\noalign{\smallskip}
				$N$	& $3134$ & $12551$ & $50184$ &  $200734$\\ 
				\noalign{\smallskip}
				\hline
				\noalign{\smallskip}
				$t_{block}$  & $14.8$ & $53.1$  & $184.5$  &   $1758.1$  \\ 
				\hline 
			\end{tabular}
		\end{center}
				\caption{CPU times (in seconds) obtained by running the block-based partition
			algorithm ($t_{block}$) }
		\label{tabe_6}
	\end{table}
	
	\section{Applications} \label{appl}
	In this section we analyze two main applications of our  partitioning 
	structure. Precisely for 2D data sets we briefly illustrate the importance of having such a versatile tool in biomathematics to assess the domains of attraction in dynamical systems.  Then for 3D data sets it will be pointed out that the flexibility with respect to the domain of  our partitioning procedure leads to an important application, i.e. modelling implicit surfaces.
	
	\subsection{Surface approximation from biomathematics}
	It is well-known that in dynamical systems saddle points partition the domain into basins of attraction of the remaining locally stable equilibria. This situation is rather common especially in population dynamics models, like competition systems \cite{C-D-P-V}. Trajectories with different initial conditions will possibly converge toward different equilibria, depending on the locations of their respective initial conditions. The set of all points that taken as initial conditions will have trajectories all tending to the same equilibrium is called the basin of attraction of that equilibrium point. 
	
	We consider the following competition model \cite{gosso12}:
	\begin{equation} \label{model3d}
	\begin{array}{ll}
	\frac{ \displaystyle  dx}{ \displaystyle  dt}=p \big(1- \frac{ \displaystyle  x}{ \displaystyle  u} \big)x-axy-bxz,  & \textrm{} \\
	\vspace{.01cm}\\
	\frac{ \displaystyle  dy}{ \displaystyle  dt}=q \big(1- \frac{ \displaystyle  y}{ \displaystyle  v} \big)y-cxy-eyz, & \textrm{} \\
	\vspace{.01cm}\\
	\frac{ \displaystyle  dz}{ \displaystyle  dt}=r \big(1- \frac{ \displaystyle  z}{ \displaystyle  w} \big)z-fxz-gyz,  & \textrm{} 
	\end{array}
	\end{equation}
	where $x$, $y$ and $z$ denote the three populations, each one competing with both the other ones
	in the same environment.
	Respectively,
	$p$, $q$ and $r$ are their growth rates, $a$, $b$, $c$, $e$, $f$ and $g$ denote the competition rates,
	$u$, $v$ and $w$ are their carrying capacities.  We assume that all parameters are nonnegative.
	
	There are eight equilibrium points; here we list only those which play a 
	role in this investigation, i.e. $E_2 = (0, v, 0)$ and $E_3 = (0, 0, w)$. 
	For further details about the study of competition  models see \cite{gosso12,Murray}.
	For suitable parameters choices the system admits multistability.
	For example, with the parameters $p = 1,q = 0.5,r = 2,a = 1, b = 2, c =0.3, e = 1,f = 3,g = 2,u = 1,v = 2,w = 1$, 
	the equilibria $E_2$ and $E_3$  are both stable equilibria and this suggests the existence of
	a separatrix surface.

\begin{figure}
		\centering
		\includegraphics[height=.26\textheight]{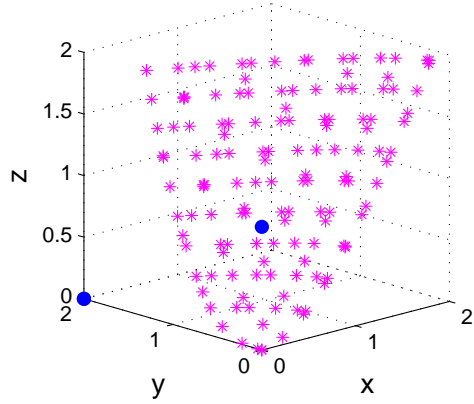} 
		\caption{A set of points lying on the surface separating $E_2$ and $E_3$ identified by blue circles.}
		\label{fig_conv}
	\end{figure}

\begin{figure}
		\centering
		\includegraphics[height=.26\textheight]{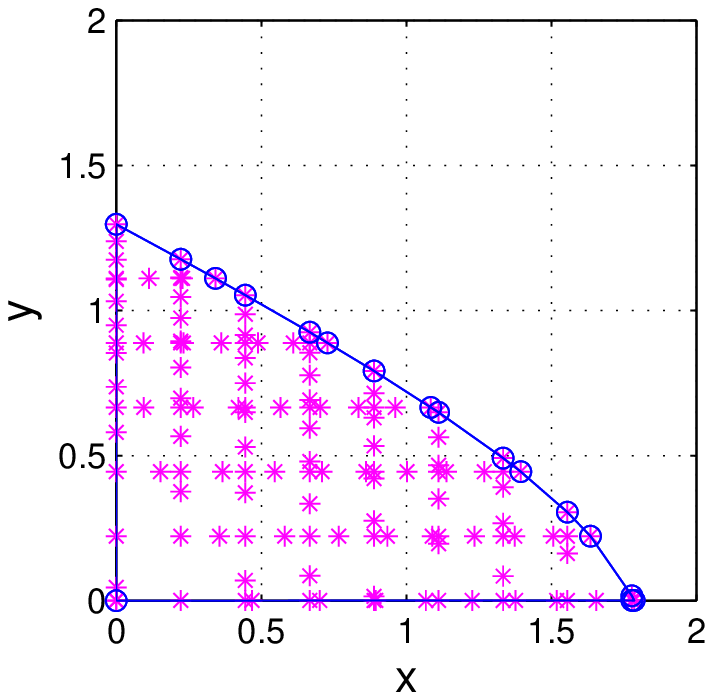} \hskip -.5cm
		\includegraphics[height=.26\textheight]{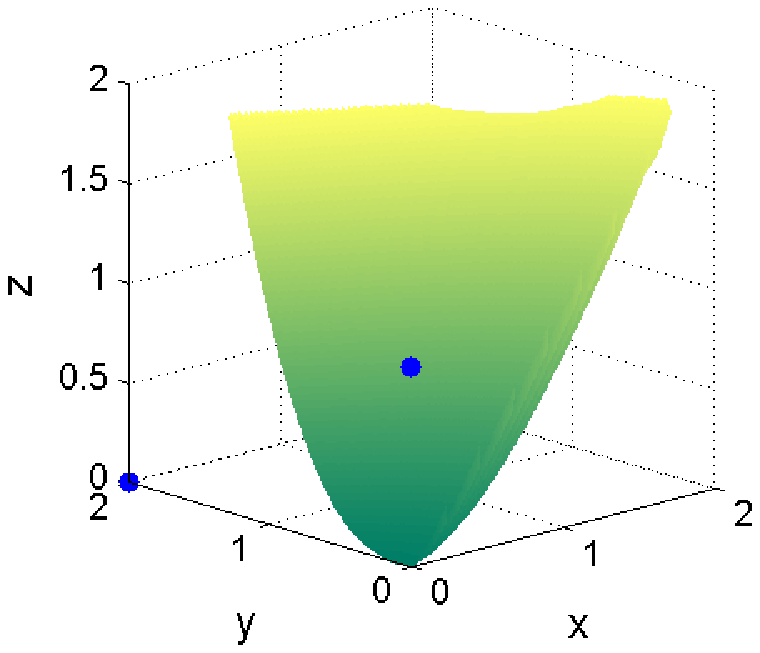}
		\caption{The convex hull defined by separatrix points (left) and the separatrix surface (right).}
		\label{fig_sep}
	\end{figure}

	To determine the separatrix surfaces for (\ref{model3d}), we need to consider a set of points as initial conditions in
	a cube domain $[0,\gamma]^3$, where $\gamma \in \RR^+$ (in the following we fix $\gamma=2$).
	Then, we take points in pairs and we check if trajectories of the two points converge to different equilibria. If this the case, we proceed with a bisection-like procedure to determine a separatrix	point. The algorithm to detect separatrix points, with a bisection routine, is analyzed in \cite{C-D-P-V}. Here we omit details and we only show points that within a certain tolerance belong to the separatrix surface in Figure \ref{fig_conv}.	Then, once the detection routine provides the separatrix points, the algorithm described in Subsection 3.1 computes the convex hull and interpolates data sites as shown in Figure \ref{fig_sep} (left and right, respectively). The surface is reconstructed by taking the shape parameter of the Wendland's function equal to $0.1$.
	
	From this sketch it is evident, especially in applications where the location of points is not in general a priori known, the importance of having a geometry-in\-de\-pen\-dent efficient partitioning structure.  	

	\subsection{Reconstruction of 3D objects}
	We present an implicit approach via partition of unity interpolation for the reconstruction of 3D objects or more in general implicit surfaces. A common problem in computer aided  design  and computer graphics is the reconstruction of surfaces defined in  terms of a set of unorganized, irregular points in 3D. Such applications also arise in computer graphics, modeling complicated 3D objects or in medical imaging \cite{Fasshauer,Wendland02}.
	
	In the approximation of 3D objects a large set ${\cal X}_N= \{ \boldsymbol{x}_i \in \mathbb{R}^{3},  i=1, \ldots, N \}$ of scattered data points, named cloud data, is generally involved. Such points should be thought as data sites belonging to an unknown two dimensional manifold $\mathscr{M}$, namely a surface in $ \mathbb{R}^{3}$.
	Thus given the set $ {\cal X}_N$,  we seek another surface $ \mathscr{M}^{*}$ that is a reasonable 
	approximation to $\mathscr{M}$. Therefore here we use an implicit approach, i.e. $\mathscr{M}$ is defined as the surface of all points $\boldsymbol{x}  \in \mathbb{R}^{3}$ satisfying the implicit equation:
	\begin{equation}
	\label{eqimp}
	f(\boldsymbol{x} )=0,
	\end{equation}
	for some function $f$, which implicitly defines the surface $\mathscr{M}$ \cite{Fasshauer,Zhu}. This means that the equation \eqref{eqimp} is the zero iso-surface of the trivariate function $f$, and therefore this iso-surface coincides with $\mathscr{M}$ \cite{Fasshauer,Wendland02}.
	The key to finding the interpolant of the trivariate function $f$,
	from the given data points $ \boldsymbol{x}_i,$ $ i=1, \ldots, N,$ is 
	to use additional significant interpolation conditions, i.e. to add an 
	extra set of \emph{off-surface points}. When the augmented data set is defined, we can then compute a three dimensional interpolant $\cal{I}$ to the total set of points \cite{Fasshauer}.
	
	In order to build the extra set of  off-surface points, 
	we assume that in addition to the point cloud data
	the set of surface oriented normals 
	$\boldsymbol{n}_i \in \mathbb{R}^{3} $ to the surface
	$ \mathscr{M}$  at the points  $ \boldsymbol{x}_i$ is also given.
	Thus we construct the extra off-surface points by
	taking a small step away along the surface normals,
	i.e. we obtain for each data point $\boldsymbol{x}_i$ two additional off-surface points, which lie \emph{outside} and \emph{inside} the manifold $ \mathscr{M}$:
	\begin{align*}
	\boldsymbol{x}_{N+i}=\boldsymbol{x}_i+ \delta \boldsymbol{n}_i, \qquad \boldsymbol{x}_{2N+i}= \boldsymbol{x}_i- \delta \boldsymbol{n}_i,
	\end{align*}
	$\delta$ being the stepsize. Note that if we have zero normals in the given normal data set, we must exclude such normals \cite{Fasshauer}. The union of the sets 
	$\cal{X}_{ \delta}^{+}=$ $  \{ \boldsymbol{x}_{N+1},\ldots, \boldsymbol{x}_{2N} \}$, $\cal{X}_{ \delta}^{-}=$ $ \{ \boldsymbol{x}_{2N+1},\ldots,$ $\boldsymbol{x}_{3N} \}$ 
	and ${\cal X}_N$  gives the overall set of points on
	which the interpolation conditions are assigned. 
	
	Now, after creating the data set, we compute the interpolant $\cal{I}$
	whose zero contour (iso-surface $\cal{I}=$ $0$) interpolates 
	the given point cloud data, and whose iso-surface
	$\cal{I}=$ $1$ and $\cal{I}=$ $-1$ interpolate 
	$\cal{X}_{ \delta}^{+}$ and $\cal{X}_{ \delta}^{-}$, respectively \cite{carr97,Fasshauer}.
	
	This problem is now reduced to a full 3D interpolation problem. Moreover
	we can notice that the large initial data set of point cloud data is 
	significantly augmented by the extra set of  off-surface points.
	Thus, from such consideration, it is evident the importance of 
	having an efficient tool which allows to compute 3D
	objects, especially in case adaptive methods are developed in the
	approximation of implicit surfaces, as in  \cite{Ohtake}.
	
	After computing the interpolant, we just render the resulting approximating surface $ \mathscr{M}^{*}$ as the zero contour of the 3D interpolant \cite{Fasshauer}. If the normals are not explicitly given, some techniques to estimate 
	the latter are illustrated in \cite{hoppe92,Wendland02}. 
	
	Here we show, several numerical experiments.
	The data sets used in our examples correspond to various point cloud data set of the Stanford Bunny.\VerbatimFootnotes\footnote{The data sets of the Stanford bunny are available at \verbdef\demo{http://graphics.stanford.edu/data/3Dscanrep/}\demo.}
	
	To approximate the 3D object, we use as local approximant the Wu's $C^4$ CSRBF that is strictly positive definite in $ \mathbb{R}^3$ cite{Wu}:
	\begin{align*}
	\phi(r)=(1-\varepsilon r)^{6}_{+} (5 \varepsilon^5r^5+30 \varepsilon^4r^4+72\varepsilon^3r^3+82\varepsilon ^2r^2 +36\varepsilon r +6). 
	\end{align*}
	
	In Table \ref{tabe_8} we report the CPU times obtained by running the block-based partition algorithm for four different sets of point cloud data.
	We remark that the interpolation conditions are almost three times larger than the original data set consisting of $N$ points (not exactly every point has a normal vector associated with it, since zero normals must be excluded).

\begin{figure}
		\begin{center}
			\includegraphics[height=.26\textheight]{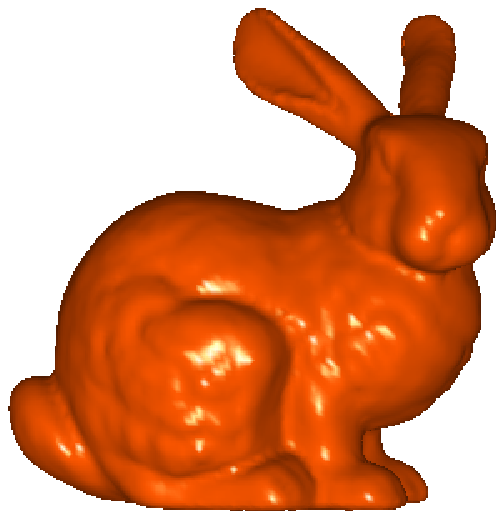} \hskip -.5cm
			\includegraphics[height=.26\textheight]{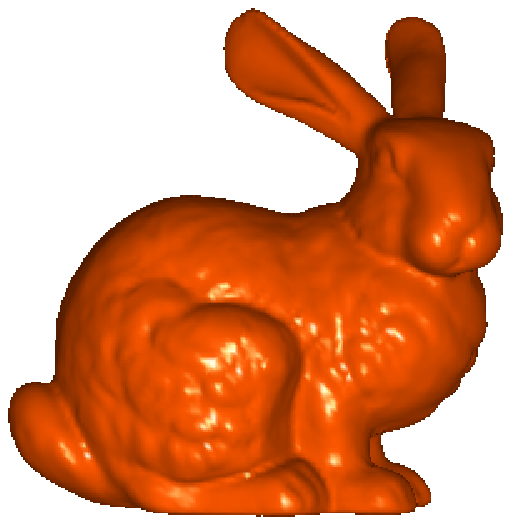} 
			\caption{The Stanford Bunny with $8171$ (left) and $35947$ (right) data points.}
			\label{bunny}
		\end{center}
	\end{figure}

	The results of the approximation algorithm for two different data sets are shown in Figure \ref{bunny}. They are obtained by taking the shape parameter $\varepsilon$ of the Wu's function equal to $0.1$ and a  grid of $100 \times 100 \times 100$  evaluation points on ${\cal R}$.\\

		\begin{table}[ht]
			\begin{center}
				\begin{tabular}{cccccc} 
					\hline\noalign{\smallskip}
					$N$	& $453$ & $1889$ & $8171$ &  $35947$ \\
					\noalign{\smallskip}
					\hline
					\noalign{\smallskip}
					$t_{block}$  & $19.4$ & $78.5$  & $642.8$  &   $7488.1$   \\
					\hline 
				\end{tabular}
			\end{center}
			\caption{CPU times (in seconds) obtained by running the block-based algorithm ($t_{block}$).}
			\label{tabe_8}
		\end{table}
	
	\section{Conclusions and work in progress}
	\label{conclusion}
	In this paper we present an efficient construction of the partition of unity interpolant. In fact the search of nearest points in the localized process, taking advantage of the block-based structure can be performed in a constant time. This is mainly due to the fact that the proposed partitioning routine is strictly related to the partition of unity subdomains, differently  from other routines.
	Extensive numerical tests and a reliable complexity analysis support our findings. 
	
	Furthermore, considering some applications in geometric modeling we show the versatility of our software.  Moreover, even if  here we deal with convex domains, it has been pointed out that the numerical tool presented can be easily adapted, with really few changes, for working in case of non convex domains, as in \cite{Heryudono}. 
	Differently from such paper, our approach excludes the employment of techniques, as for example conformal maps. In fact in that paper the Schwarz-Christoffel transformation is used to map the interpolant, which is built exclusively in the unit disk and not on the  irregular domain, onto a known polygon. Then a global  method is considered and consequently the solution of the interpolation of large scattered data sets cannot be achieved. On the contrary, using our technique we could get an accurate solution in case of large data sets in an irregular domain and the problem solved in a relatively small time. In that case, applications of our algorithm could arise in the context of interpolation and PDEs \cite{Chen14,Deparis14,Safdari,Shcherbakov}, allowing to deal with larger data sets in a reasonable time. Further investigations in this direction are needed.

	Moreover, work in progress also consists in extending the proposed block-based partitioning scheme so that it allows to consider subdomains having variable radii. This turns out to be meaningful especially when strongly non-uniform data are considered. 

	\section*{Acknowledgements}		
		The authors acknowledge financial support from the GNCS--INdAM.
		
		

\end{document}